\documentclass{article}

\usepackage{amssymb,amsmath,amsthm}
\usepackage{graphicx}
\usepackage[english]{babel}
\usepackage{enumitem}
\usepackage{color}
\usepackage{a4wide}

\def\rr{\mathbb{R}} 
\def\sF{\mathcal{F}} 
\def\Pr{\mathbb{P}} 
\def\dif{\mathrm{d}} 
\def\landau{\mathcal{O}} 
\def\part{\mathcal{T}_{\mathbf{h}}^N} 
\def\e{\textnormal{e}} 

\def\ntilde{\widetilde}
\def\nhat{\widehat}
\def\nbar{\overline}

\newtheorem{algorithm}{Algorithm}
\newtheorem{definition}{Definition}
\newtheorem{theorem}{Theorem}
\newtheorem{corollary}{Corollary}

\begin{document}

\renewcommand{\topfraction}{0.85}
\renewcommand{\textfraction}{0.1}
\renewcommand{\floatpagefraction}{0.75}

\parindent 0pt

\title{Almost sure convergence of numerical approximations for Piecewise Deterministic Markov Processes}

\author{%
{\sc
Martin G.~Riedler\thanks{Email: martin.riedler@jku.at}} \\[2pt]
Institute for Stochastics, Johannes Kepler University,\\[2pt] Altenbergerstra{\ss}e 69, 4040 Linz, Austria
}

\maketitle

\begin{abstract}
{Hybrid systems, and Piecewise Deterministic Markov Processes in particular, are widely used to model and numerically study systems exhibiting multiple time scales in biochemical reaction kinetics and related areas. In this paper an almost sure convergence analysis for numerical simulation algorithms for Piecewise Deterministic Markov Processes is presented. The discussed numerical methods arise through discretisina a constructive method defining these processes. The stochastic problem of simulating the random, path-dependent jump times of such processes is reformulated as a hitting time problem for a system of ordinary differential equations with random threshold. Then deterministic continuous methods (methods with dense output) are serially employed to solve these problems numerically. We show that the almost sure asymptotic convergence rate of the stochastic algorithm is identical to the order of the embedded deterministic method. We illustrate our theoretical findings by numerical examples from mathematical neuroscience were Piecewise Deterministic Markov Processes are used as biophysically accurate stochastic models of neuronal membranes.}
{stochastic simulation; hybrid algorithm; almost sure convergence; Piecewise Deterministic Markov Process.}
\end{abstract}


\section{Introduction}
%
%
%
%
%
In recent years the number of applications of hybrid stochastic processes to model systems in biology, (bio-)chemical reaction kinetics and mathematical neuroscience have increased rapidly. These models either stem from a direct modelling approach, e.g., models of excitable biological membranes \cite{BuckwarRiedler,ClayDefelice,Wainrib1,RudigerFalcke}, or result from a multiscale approximation to more accurate particle models exhibiting clearly separated time scales, e.g., in biochemical reaction systems \cite{Huisinga,Haseltine,Kalantzis,SalisKaznessis} or in gene regulatory networks \cite{Zeiseretal}. The mathematically correct treatment of these is within the framework of \emph{General Stochastic Hybrid Systems} \cite{Bujorianu}. Additional to these novel `bioscience' applications for hybrid systems there are also more `classical' applications of high interest in fields such as control and queueing theory, or models in financial mathematics and ecology, cf.~\cite{Bujorianu,Davis1,Davis2,Jacobsen,Hanson} and references therein.
%
%
%
%
%
In this paper we focus on  Piecewise Deterministic Markov Processes (PDMPs), which are an important class of hybrid systems including most of the hybrid models considered in the literature, see \cite{Huisinga,BuckwarRiedler,ClayDefelice,Davis1,Davis2,Haseltine,Jacobsen,Kalantzis,Wainrib1,RudigerFalcke,Hanson,Zeiseretal} out of the above mentioned studies. 
PDMPs are strong Markov processes that combine continuous deterministic time evolution and discontinuous, instantaneous, random `jump' events. Specifically, the dynamics of the two components are intrinsically intertwined. On the one hand, the continuous time-evolution, defined by ordinary differential equations (ODEs), depends on the outcomes of discrete events via randomly changing parameters and, on the other hand, the probability of the discrete events happening, i.e., a random, instantaneous change in a parameter, depends on the time-evolution -- the path -- of the continuous variables.

%
%
%
%
Due to the high complexity of hybrid models, particularly in biochemical applications, they are studied extensively and almost exclusively by numerical means. To this end either so called `pseudo-exact' algorithms, i.e., simple statistically approximate algorithms obtained by an ad-hoc model approximation, or statistically exact algorithms are employed. Here \emph{statistically exact} means that in principle the simulation method produces paths that are samples of the distribution of the underlying stochastic process in contrast to pseudo-exact methods which produce sample paths from a distribution approximating the distribution of the underlying stochastic process. However, even statistically exact algorithms are \emph{exact only in theory due to a numerical error which arises inevitably} in solving most systems of differential equations in actual implementations. So far model accuracy was primarily of interested as the main source of error and thus considerations regarding the numerical error are in general neglected, cf., e.g., \cite{Huisinga,Haseltine,Kalantzis}. Yet ultimately, even if a theoretically exact PDMP formulation of the model is considered or the model is highly accurate, \emph{numerical studies are conducted by numerical approximations} to the PDMPs as an analytic representation of the paths is in general not available. Despite the importance and widespread use of numerical studies a numerical analysis, in particular, a thorough analytical investigation of the convergence and error behaviour of algorithms to be used, is still missing. The aim of this article is to provide a -- to the best of our knowledge -- first contribution towards this goal. We note that the statistically exact algorithms that were introduced in the studies \cite{Huisinga} and \cite{Zeiseretal} fall within the class we consider, however, in the present study a convergence analysis is carried out.

%
%
%
%
%
%
%
In particular in the present study we are interested in the convergence of numerical approximations to PDMPs in a pathwise sense which corresponds to the fact that numerical simulations are carried out path by path. 
The methods we present for approximating a PDMP incorporate as an integral part \emph{continuous ODE methods}, also called methods with dense output.
Apart from numerically solving the deterministic inter-jump dynamics the key problem in simulating a PDMP is simulating the random, path-dependent jump times. As shown in Section \ref{section_convthm} this problem can be reformulated and then, combined with the numerical solution of the inter-jump dynamics, yields a \emph{hitting time problem with random threshold}. This we solve using continuous ODE methods. The main feature of such a method is that it does not only provide a numerical approximation to the exact solution at discrete grid points but provides an approximation of the whole path over the whole interval. Essentially, continuous methods are approximations on a discrete grid with an interpolation formula for the intervals between the grid points. 
Hence, these methods are naturally suited for solving hitting time problems.
The main result of this paper is that numerical approximations of PDMPs built on continuous ODE methods conserve the order of convergence of the underlying continuous method. That is, if an approximation is constructed using, e.g., a continuous Runge-Kutta method of order $p$, then also the almost sure convergence of the stochastic approximation to the PDMP is of order $p$.

As for the PDMPs discussed in the present study we distinguish several types. The main part of the paper is devoted to processes with jumps occurring only in a fixed subset of the components of a vector-valued PDMP which are otherwise piecewise constant. Furthermore we assume that the jump heights are discretely distributed. We consider this particular structure for processes as these arise in applications in mathematical neuroscience, which initially motivated this study, and as a multiscale approximation to certain chemical reaction systems. In these models fast and slow reactions, modelled with reaction rate equations and instantaneous random jumps, respectively, affect different set of reactants, however, each possesses rates also depending on the other set of reactants. As jumps corresponds to an instantaneous event where one or more individual reactants change their state, in these models jump heights are typically integer valued. The second class of PDMPs comprises processes for which jumps may occur in all components and these components need not be necessarily piecewise constant, but jumps are still discretely distributed. These PDMPs include chemical reaction systems where fast and slow reactions may affect the same type of reactants. Finally, we also consider this last general type of processes with continuously distributed jump heights.
\medskip

%
%
%
%
The remainder of the paper is organised as follows. In Section \ref{section_PDMP} we present a brief definition of PDMPs and their construction from sequences of independent, identically distributed (i.i.d.) standard uniform random variables. This provides a theoretically exact simulation algorithm which takes the role of an exact solution for our convergence analysis. Section \ref{section_convthm} presents the approximate algorithms we consider in this paper and also contains the main convergence theorem. The proof of the convergence theorem can be found in Section \ref{section_proof}. We extend the convergence result to other classes of PDMPs in Section \ref{section_extension}. Finally, in Section \ref{section_examples} we present some numerical experiments using examples from the neuroscience literature to illustrate the theoretical findings and draw some conclusions for the usage and implementations of the numerical methods.

\section{PDMPs and their construction from i.i.d.~sequences of uniform random variables}\label{section_PDMP}

In this section we give a brief introduction to PDMPs, introducing the objects used for their definition and the construction of their paths. In particular, we present Algorithm \ref{exact_alg} which provides the `exact solution' the numerical approximations we present converge to. For a general discussion of PDMPs we refer to the monographs \cite{Davis2} and, more recently, \cite{Jacobsen} or the PhD thesis of the present author \cite{RiedlerPhD}.
\medskip

For the main part of the paper we consider a PDMP $(X(t))_{t\in [0,T]}$ $=(X(t,\omega))_{t\in[0,T]}$ to consist of two qualitatively different components, i.e., $X(t)=(Y(t),\theta(t))\in\rr^{d+m}$. For these we assume that $Y(t)$ possesses continuous paths in a set $D\subseteq\rr^d$ and $\theta(t)$ is right continuous and piecewise constant in $K\subset\rr^m$, where $K$ is an at most countable set. We denote by $t_n$, $n\geq 1$, the jump times of the component $\theta(t)$ and by $E=D\times K$ the state space of the process. Such a PDMP is uniquely defined (up to versions) by its \emph{characteristic triple} $(f,\lambda,\mu)$, where the component $f$ is used to define the deterministic, continuous evolution of the trajectories in between jumps and the components $\lambda$ and $\mu$ yield the stochastic jump mechanism. In detail these objects are: The first is a measurable map $f:E\to\rr^d$ such that for each $\theta\in K$ the differential equation
\begin{equation}\label{ode_set}
\dot{y}=f(y,\theta)
\end{equation}
possesses a unique solution on $[0,T]$ with values in $D$ for all initial conditions $y(0)=y_0\in D$. For notational convenience we introduce the map $t\mapsto\phi(t,x_0)$ which denotes the unique solution to the system
\begin{equation}\label{IVP}
 \left(\!\!\! \begin{array}{c} \dot{y}\\ \dot{\theta}\end{array}\!\!\!\right)=\left(\!\!\!\begin{array}{c} f(y,\theta)\\ 0\end{array}\!\!\!\right)
\end{equation}
with respect to the initial value $x_0=(y_0,\theta_0)\in E$. It satisfies $\phi(t,x_0)\in E$ for all $t\in[0,T]$ and all initial values $x_0\in E$. The system \eqref{IVP} defines the paths of the PDMP in between jumps, see Algorithm \ref{exact_alg} below. Thus we can think of such a PDMP as an ODE system with a parameter changing randomly in a certain way now described: The non-negative, measurable map $\lambda:E\to \rr_+$ is the instantaneous intensity of the jumps. It defines a probability distribution $\rr_+$ for all $x\in E$ via the survivor function\footnote{A survivor function $t\mapsto S(t)\in[0,1]$, $t\geq 0$, of a non-negative random variable states the probability that this random variable is larger than $t$, for discussion in relation to PDMPs see \cite[Sec.~4.1]{Jacobsen}.}
\begin{eqnarray}\label{survfunc}
S(t,x\bigr)&\,=\,&\exp\biggl(-\int_0^t \lambda\bigl(\phi(s,x)\bigr)\,\dif s\biggr)\,,
\end{eqnarray}
which is used to define the distribution of the inter-jump times. That is, $S(t,x)$ states the probability that there does not occur a jump in the next time interval of length $t$ conditional on the process being currently in state $x$. We always assume that $\lambda$ is path-integrable, i.e.,
\begin{equation*}
 \int_0^T \lambda(\phi(s,x))\,\dif s<\infty\quad\forall\,x\in E\,.
\end{equation*}
Finally, $\mu$ is a Markov kernel from $E$ into $\rr^m$ defining the distribution of the jump heights in the component $\theta(t)$ conditional on the pre-jump values. We assume that $\mu(x,\cdot)$ is a discrete probability measure satisfying $\mu\bigl((y,\theta),\{0\}\bigr)=0$ for all $x=(y,\theta)\in E$ and that $\mu$ is continuous in $x$, i.e., $x\mapsto \mu(x,A)$ is continuous for every Borel set $A\subset\rr^m$. Clearly, the at most countable support $K_x$ of the measure $\mu(x,\cdot)$ has to be such that the resulting post-jump values are in $K$, i.e., for all $x=(y,\theta)$ it holds that $\theta+K_x\subseteq K$.\medskip

Algorithm \ref{exact_alg} below can be used to prove the existence of a unique PDMP to given characteristics in a constructive way \cite{Davis2}. More importantly, this algorithm also provides a theoretically exact method for simulating paths of a PDMP on a probability space that supports a sequence of i.i.d.~standard uniform random variables. Here `exact' denotes the fact that the distribution of the process defined by Algorithm \ref{exact_alg} equals the distribution of the PDMP to the triple $(f,\lambda,\mu)$. However, the adverb `theoretically' should emphasise the fact that in general the algorithm is an `impossible' algorithm as neither the initial value problem (IVP) \eqref{IVP} nor the implicit equation \eqref{next_jump_time} arising in Algorithm \ref{exact_alg} can usually be solved exactly. Hence, the algorithm cannot be employed in practice to simulate trajectories. Yet, it can be used as a comparison to approximate algorithms and thus, in typical terms of numerical analysis, plays the role of an exact solution for a convergence analysis. Now, let $(\Omega,\sF,\Pr)$ be a probability space which supports a sequence of i.i.d.~standard uniform random variables $U_n$, $n\in\mathbb{N}$.

\begin{algorithm}\label{exact_alg} An exact simulation algorithm for a PDMP is given by:
\setlength{\leftmargini}{3.5em}
\begin{enumerate}[align=left]
\setlength{\leftmargini}{3.5em}
  \setlength{\labelwidth}{2.5em}
  \setlength{\labelsep}{1.0em}
 \item[Step 1.]  Fix the initial time $t_0=0$ and initial condition \mbox{$X(0)=x_0=(y_0,\theta_0)\in E$} and set a jump counter $n=1$.
 \item[Step 2.] Simulate a uniformly distributed random variable $U_{2n-1}$ and solve
\begin{equation}\label{next_jump_time}
 S(\tau,X(t_n))=U_{2n-1}
\end{equation}
with respect to $\tau$ to obtain the waiting time until the next jump time, i.e., $\tau=t_{n}-t_{n-1}$. Then for $t\in(t_{n-1},t_n]$ set
\begin{equation*}
 X(t)=\phi(t-t_{n-1},X(t_{n-1}))\,. 
\end{equation*}
If $t_{n}\geq T$, stop at $t=T$.
\item[Step 3.] Otherwise, simulate a post-jump value $\theta_{n}$ for the piecewise constant component according to the distribution of the jump heights given by $\mu(\phi(t_{n}-t_{n-1},X(t_{n-1})),\,\cdot\,)$ via the uniformly distributed random variable $U_{2n}$ and set
\begin{equation*} X(t_{n})=\left(\!\!\begin{array}{c} \phi(t_{n}-t_{n-1},X(t_{n-1})) \\ \theta_{n}\end{array}\!\!\right).\end{equation*}
 \item[Step 4.] Set $n=n+1$ and start again with Step 2.
\end{enumerate}
\end{algorithm}

Algorithm \ref{exact_alg} is almost surely well-defined unless one of the following two cases occur: first, if neither a unique nor a solution at all exists to equation \eqref{next_jump_time}, which can occur for PDMPs, and second, if trajectories possess countably many jumps before time $T$ with positive probability. In this study the first case is excluded by the conditions on $\lambda$ in the convergence theorem and hence for the remainder of the paper we can assume that a unique solution to \eqref{next_jump_time} always exists. The second case is excluded if we assume that the characteristics $(f,\lambda,\mu)$ are such that the PDMP is a \emph{regular} jump process. That is, there are almost surely only finitely many jumps in any finite time interval. A simple condition which guarantees the regularity of a PDMP is that the jump rate $\lambda$ is bounded, see \cite{Davis2,Jacobsen} which also contain more general conditions.

We remark that the class of PDMPs to triples $(f,\lambda,\mu)$ as given above contains some prominent classes of stochastic processes as special cases. 
A first special case is given by ODEs with Markovian switching. For these processes the jump rate is given by jumps of a Poisson process that is independent of the paths followed in between jumps, i.e., $\lambda\equiv const$. Such equations can be written as stochastic differential equations driven by a Poisson process, i.e., stochastic differential equations of jump type (JSDEs) without a diffusion term, cf.~\cite{MaoYuan}.
In general, JSDEs without a diffusion part are special cases of PDMPs \cite{Jacobsen}.
Secondly, PDMPs are generalisations of continuous-time Markov chains.
These are piecewise constant processes and thus follow a trivial evolution between successive jump times, i.e., $f\equiv 0$.
However, in contrast to switching ODEs the jump intensity may be state dependent, i.e., $\lambda\not\equiv const$. Such processes are most importantly used in stochastic models of chemical reaction kinetics and are treated numerically by the Stochastic Simulation Algorithm (SSA) \cite{Gillespie}. %
%
%
%
%
We note that the numerical advantage of these special cases is that inter-jump times can be simulated exactly by sampling from an exponential distribution, that is, the implicit equation \eqref{next_jump_time} can be solved exactly. Accordingly there exists a large amount of literature on simulation methods for these processes and a numerical analysis thereof, see, for instance, references in the recent studies \cite{BuckwarRiedlerJSDEs,SSASummary}.\medskip

\textbf{Example.} To conclude this brief introduction we present a concrete example for a PDMP arising in mathematical neuroscience which is describe in more detail in Section \ref{section_examples}. The component $Y(t)$ is one-dimensional taking values in $[0,E_{\textnormal{Na}}]$ for initial conditions therein and the component $\theta(t)$ is 8--dimensional where $K=\{\theta\in\{0,\ldots,N\}^8\,:\ \sum_{k=1}^8 \theta_k=N\}$ with $N\in\mathbb{N}$. Thus the phase space is $E=[0,E_{\textnormal{Na}}]\times K$ and the characteristics are given as follows. The family of ODEs \eqref{ode_set} is given by
\begin{equation*}
C\dot y\,=\, -\nbar g_{\textnormal{Na}}\,\theta_{8}\,(y-E_{\textnormal{Na}})-\nbar g_\textnormal{L} y\,,\quad \theta\in K
\end{equation*}
with constants $\nbar g_{\textnormal{Na}},\, E_{\textnormal{Na}},\, \nbar g_\textnormal{L} > 0$, and the stochastic dynamics are given by the jump-rate
\begin{equation*}
\lambda\bigl((y,\theta)\bigr)=
\left(\!\!\begin{array}{c}
a_m(y) \\ b_m(y) \\ a_h(y) \\ b_h(y)
      \end{array}\!\!\right)^T\,
\left(\begin{array}{cccccccc}
 3 & 2 & 1 & 0 & 3 & 2 & 1 & 0\\
 0 & 1 & 2 & 3 & 0 & 1 & 2 & 3\\
 1 & 1 & 1 & 1 & 0 & 0 & 0 & 0\\
 0 & 0 & 0 & 0 & 1 & 1 & 1 & 1\\
\end{array}\right)\,
\left(\!\!\begin{array}{c}
\theta_1 \\ \vdots \\ \theta_8
      \end{array}\!\!\right)
\,.
\end{equation*}
and the transition measure $\mu$ consisting of point probabilities of the form
\begin{equation*} 
\mu\bigl((y,\theta),\{(-1,1,0,0,0,0,0,0)^T\}\bigr)=\frac{3\,a_m(y)\,\theta_1}{\lambda\bigl((y,\theta)\bigr)}\,.
\end{equation*}

\subsection{A property of the jump height distribution}\label{section_simulations}

\begin{figure}
\centering
\includegraphics[width=0.7\textwidth, clip=true, trim=80 80 100 50]{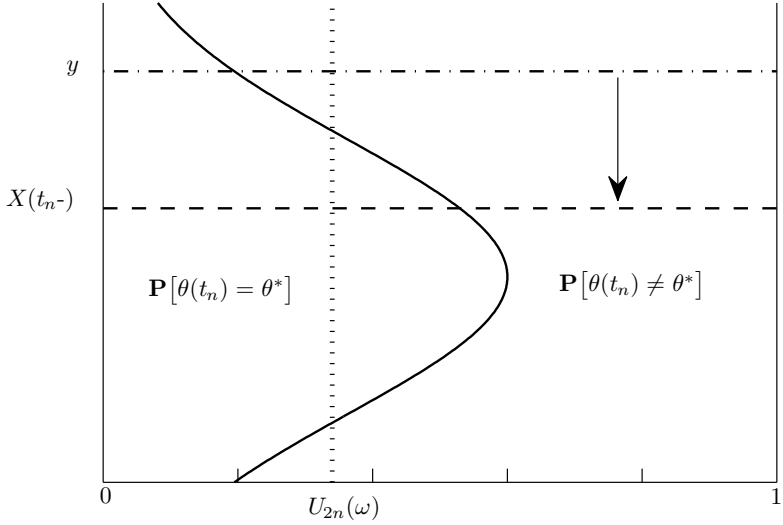}
\caption{Schematic view of the domain of $H$ for probabilities $\mu(X(t_n-),\cdot)$ with discrete support. The curve marks points of discontinuity. The function $H$ is constant $\theta^\ast$ for arguments $(x,u)$ left of the curve and piecewise constant with different values right of it. $\theta^\ast$ is the realised $n$th jump height and we see $H(y,U_{2n}(\omega))=\theta^\ast$ for $y$ close enough to $X(t_n-)$ as the probability that $U_{2n}$ lies on a discontinuity is zero.}\label{jump_transitions}
\end{figure}

In this section we derive a consequence of the continuity of $x\mapsto\mu(x,\cdot)$ which is important for the proof of the convergence result. First, we recall that a standard method used to define random variables from uniform distributed random variables as is needed in Algorithm \ref{exact_alg} is the \emph{inverse CDF method}, see, e.g., \cite{Gentle}. In particular, solving \eqref{next_jump_time} in Step 2 of Algorithm \ref{exact_alg} is a version of the inverse CDF method for a distribution on the positive reals. We assume that also for the simulation of the jump heights a version of the inverse CDF method is employed, that is, there exists a function $H:E\times(0,1)\to K$ such that for all measurable $A\subset\rr^m$ it holds
\begin{equation*}
\Pr\bigl[H(x,U)\in A\bigr]=\mu(x,A)\,, 
\end{equation*}
cf., e.g., \cite[Corol.~23.4]{Davis2}. Hence, the post jump value $\theta_{n+1}$ in Step 3 of Algorithm \ref{exact_alg} is given by $\theta_n+H(\phi(t_{n+1}-t_n,X(t_n)),U_{2n+2})$. As $K$ is an at most countable set  there exists another at most countable set $K'$ containing all supports $K_x$ of the distributions $\mu(x,\cdot)$ for all $x\in E$. Thus the Markov kernel $\mu$ can be written as a linear combination of Dirac measures
\begin{equation*}
\mu(x,\cdot)=\sum_{k\in K'} p_k(x)\delta_k, 
\end{equation*}
where $p_k(x)$ is the point probability of the jump height $k\in K'$ conditional on the realised pre-jump value $x$. The probabilities $p_k(x)$ may very well be -- and in general are -- zero for certain pairs $(x,k)$. We index the events in $K'$ from $1$ to $|K'|\leq \infty$ and obtain a function $H$ by defining
\begin{equation}\label{definition_of_H}
H(x,u)=k_i \ \textnormal{ if } \sum_{j=1}^{i-1} p_{k_j}(x) \leq u < \sum_{j=1}^i p_{k_j}(x)\,. 
\end{equation}
The definition $H$ as in \eqref{definition_of_H} refers to the standard simulation method for discrete probabilities from standard uniform random variables by `binning' the interval $(0,1)$, see, e.g., \cite{Gentle}. However, the property that $x\mapsto \mu(x,A)$ is continuous for all measurable sets $A$ now implies that the point probabilities $p_k(x)$ depend continuously on $x$. Thus $H$ is continuous in the first component at a point $x$, i.e., $\lim_{z\to x}H(z,u)=H(x,u)$, unless $u=\sum_{j=1}^i p_{k_i}(x)$ for some $i=1,\ldots,|K'|$. That is, unless $u$ equals a discontinuity point of $u\mapsto H(x,u)$. However, the set of discontinuity points is at most countable and hence the probability of a standard uniform random variable $U$ taking such a value is zero. Therefore it holds that we can choose almost surely a $z$ close enough to $x$, say $|z-x|\leq\kappa$, where $\kappa=\kappa(x,U)$ depends on $x$ and $U$, such that
\begin{equation}\label{continuity_of_H}
|H(x,U)-H(z,U)|=0\,, 
\end{equation}
i.e., for $z$ close enough to $x$ the realised jump heights coincide, cf.~Fig \ref{jump_transitions}.

Finally, we remark, that the function $H$ is in general not unique. However, choosing a different $H$ just results in a different version of the same PDMP on the same probability space. Thus throughout the paper we assume that one such $H$ is fixed which is used in Step 3 of Algorithm \ref{exact_alg} and also in Step 3 of the numerical approximation Algorithm \ref{approx_alg} in Section \ref{section_convthm} and it satisfies \eqref{continuity_of_H}.

\section{The main convergence theorem}\label{section_convthm}
In this section we first precisely state and discuss the convergence concept we are interested in before we continue with constructing the numerical methods by discretising Algorithm \ref{exact_alg}. Subsequently, at the end of the section we state the central result of this study: the almost sure convergence of the numerical approximations and their asymptotic order of convergence. \medskip

The set-up is as in Section \ref{section_PDMP}, i.e., $(\Omega,\sF,\Pr)$ is a probability space supporting a sequence of i.i.d.~standard uniform random variables  $U_n,\,n\geq 1$. Further, we fix a finite time interval $[0,T]$. Then for given characteristics $(f,\lambda,\mu)$ as in Section \ref{section_PDMP} we use the sequence of uniform random variables to construct a regular PDMP $(X(t,\omega))_{t\in[0,T]}$ on this probability space via Algorithm \ref{exact_alg}.  
Moreover, we assume that numerical approximations $(\nhat X(t,h,\omega))_{t\in[0,T]}$ are defined on the same probability space. Here $h$ denotes a defining parameter of the numerical approximation, such as a discretisation step size.
In the following we refer to the PDMP $(X(t,\omega))_{t\in[0,T]}$ as the \emph{exact PDMP} and to its numerical approximation $(\nhat X(t,h,\omega))_{t\in[0,T]}$ as the \emph{approximate PDMP}.
Further, $t_n(\omega)$ and $\nhat t_n(h,\omega)$, $n\geq 1$, denote the jump times of the exact PDMP and its approximation, respectively. Finally, the number of the exact PDMP's jumps in the time interval $[0,T]$ is denoted by $N(T,\omega)$, where $N(T,\omega)<\infty$ almost surely.

\begin{figure}
\centering
\includegraphics[width=1\textwidth, clip=true, trim=30 65 40 70]{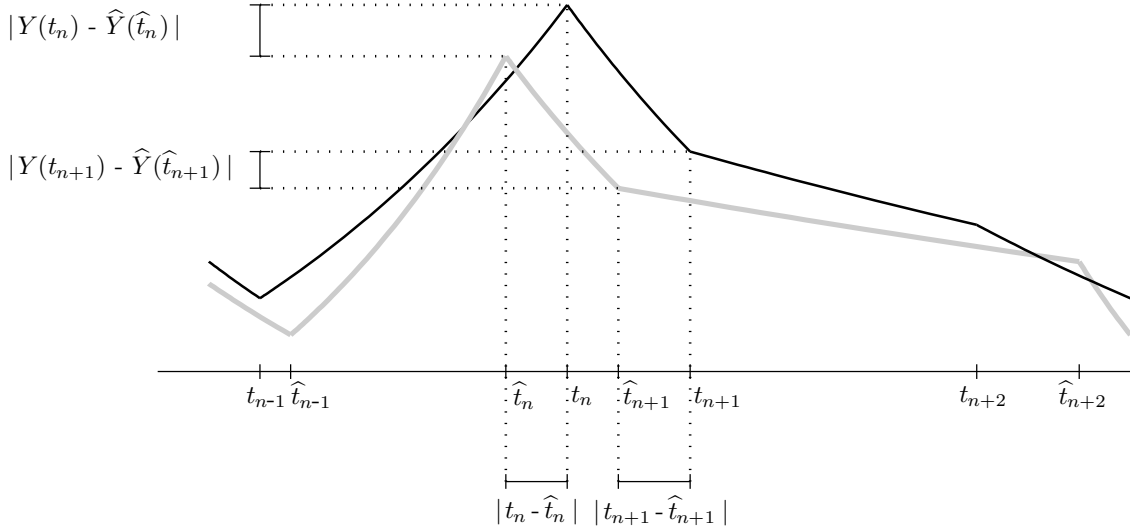}
\caption{Schematic relation of the exact PDMP's continuous component $Y(t)$ (black line) to its approximation $\nhat Y(t)$ (grey line) in phase space and the approximation errors we need to control to obtain path-wise convergence of Algorithm \ref{approx_alg}.}\label{relation}
\end{figure}

%
%
%
%
In this study we are interested in the \emph{pathwise convergence of global errors} of the form 
\begin{equation}\label{max_phase_conv}
\lim_{h\to 0} \max_{n=1,\dots, N(T,\omega)} |X(t_n(\omega),\omega)-\nhat X(\nhat t_n(h,\omega),h,\omega)| = 0
\end{equation}
and
\begin{equation}\label{endpoint_error}
\lim_{h\to 0}|X(T,\omega)-\nhat X(T,h,\omega)| = 0
\end{equation}
for almost all $\omega\in\Omega$.
That is, \eqref{max_phase_conv} and \eqref{endpoint_error} state the almost sure convergence of the approximate PDMP to the exact PDMP at their jump times and at the endpoint, cf.~Fig.~\ref{relation}.
%
%
%
%
%
It is a standard approach in numerical analysis, whether deterministic or stochastic, to define the error of a numerical approximation by a measure of the difference of the exact and approximate process at discrete time points in the approximation interval $[0,T]$ \cite{hairer,KloedenPlaten,MilsteinTretyakov}.
Here the jump times of the PDMP and its approximation are used as a discretisation of the time interval. 
%
%
%
%
However, in general, the exact and approximate jump times do differ and hence yield different time discretisations for the exact and approximate PDMP. This is in contrast to error definitions typically used in numerics.
Moreover, typically, the number and locations of grid points in error definitions are the same for each trajectory and increase in number for decreasing step size $h$. It is exactly the opposite in the error definition \eqref{max_phase_conv} as $N(T,\omega)$ varies with $\omega\in\Omega$ but is fixed over variations in $h$.

As the numerical approximation is continuous over time almost everywhere we are in principle able to use a `typical' error definition.
However, for ordinary or stochastic differential equations the time points the error is evaluated at arise naturally as the grid points of the discretisation used in the numerical method. For PDMP approximations as considered in this study this is not the case.
In general, the methods we consider do not possess the same grid points for any trajectory of one PDMP to approximate and they possibly do not even possess the same number of grid points.
Therefore, setting particular but arbitrary points apart in using them in an error definition without a mathematical reasoning does not seem preferable to us.
However, for a PDMP all pathwise numerical methods necessarily need to approximate the jump times. Therefore, relating the values of the exact and approximate PDMP at their respective jump times is reasonable insofar as an accurate approximation is expected to behave after the $n$th jump just as the exact PDMP behaves after the $n$th jump. 
%
%
%

%
%
%
%
We have already mentioned that, in general, the jump times of the exact and approximate PDMP differ. Therefore
it is not sufficient to consider the errors \eqref{max_phase_conv} and \eqref{endpoint_error} alone. They may become arbitrarily small with the exact and approximate trajectories still substantially differing as possibly only the values at the jump times converge but the jump times themselves remain separated.
Thus, we additionally require from a convergent method that also the jump times of the approximate PDMP converge to the jump times of the exact PDMP, that is, for almost all $\omega\in\Omega$
\begin{equation}\label{max_time_conv}
\lim_{h\to 0}\max_{n=1,\dots,N(T,\omega)}|t_n(\omega)-\nhat t_n(h,\omega)| =0\,.
\end{equation}

\begin{definition}\label{def_convergence} We say an approximation is pathwise convergent if \eqref{max_phase_conv}--\eqref{max_time_conv} hold. In addition, an approximation is of order $p$ if $p$ is the largest integer such that the asymptotic behaviour of \eqref{max_phase_conv}--\eqref{max_time_conv} is $\landau(h^p)$.
\end{definition}

In order to construct an approximation algorithm for the numerical simulation of a PDMP's paths we start from the theoretically exact Algorithm \ref{exact_alg}. We obtain an approximate algorithm by discretising the problems \eqref{IVP} and \eqref{next_jump_time} and solving them numerically.
A numerical solution of \eqref{IVP} can be obtained by a standard ODE method, thus substituting Step 3 with an approximation is straightforward. However, the delicate part is numerically solving \eqref{next_jump_time} for which obviously a numerical solution of \eqref{IVP} is needed as an integrand which is integrated until an \emph{a priori unspecified time} $\tau$, cf.~the definition of the survivor function \eqref{survfunc}. Hence, the two problems have to be solved in parallel including a mechanism for detecting the time $\tau$ conditional on a realisation of a standard uniform random variable $U$.
For an efficient solution we transform the equation for the next jump time \eqref{next_jump_time} into an equivalent problem. This can then be combined with the IVP \eqref{IVP} yielding a merged formulation of the two problems which allows for a numerical solution by continuous ODEs methods. The resulting numerical approximation is finally given by Algorithm \ref{approx_alg} below.\medskip

Taking the logarithm of equation \eqref{next_jump_time} we obtain the equivalent equation
\begin{equation}\label{definition_of_w}
  w(\tau,x)=\int_0^\tau\lambda(\phi(s,x))\dif s=-\log U.
\end{equation}
Thus $w(\tau,x)=-\log S(\tau,x)$ is the logarithm of the survivor function $S$. Differentiation of $w$ with respect to $\tau$ yields
\begin{equation*}
\dot w(t,x)= \lambda(\phi(t,x))\,.
\end{equation*}
This in turn yields that calculating the next jump time from \eqref{next_jump_time}, i.e., Step 2 in Algorithm \ref{exact_alg}, is equivalent to solving the IVP
\begin{equation}\label{sys2}
 \left(\!\!\! \begin{array}{c} \dot{y}\\ \dot{\theta} \\ \dot{w}\end{array}\!\!\!\right)=\left(\!\!\!\begin{array}{c} f(y,\theta)\\ 0 \\ \lambda(y,\theta)\end{array}\!\!\!\right),\quad \left(\!\!\!\begin{array}{c} y(0)\\ \theta(0)\\ w(0)\end{array}\!\!\!\right)=\left(\!\!\begin{array}{c} X(t_n) \\ 0 \end{array}\!\!\right),\ t\in[0,\infty)
\end{equation}
which is integrated until the component $w(t)$ hits the threshold $-\log U_{2n-1}$.
The hitting time
\begin{equation*}
\tau=\inf\,\{t>0:w(t)=-\log U_{2n-1}\}
\end{equation*}
is equivalent to the waiting time until the next jump, where, in general, $\tau=\infty$ may be possible unless, e.g., $\lambda$ is bounded away from zero.

\medskip

We need to remark on two aspects at this point. Firstly, the set-up of the hitting time problem \eqref{sys2} is as in \cite{Huisinga}. The authors therein state how it can be solved using continuous methods and also present an ad-hoc implementation of an event detection for this specific problem. However, as the authors are primarily focused on the modelling and simulation results, they suppose that \eqref{sys2} can be computed ``up to any desired accuracy and neglect the discretisation error'' \cite[p.~6]{Huisinga}. This may be a reasonable assumption for the solution of any standard ODE without jumps in the parameters -- as numerical methods are well studied in this case --, however, \emph{additional assumptions are needed and have to be considered} such that an event detection is possible with arbitrary accuracy, cf.~\cite{Shampine}. This aspect is not discussed by the authors. An algorithm for simulating a PDMP's path repeatedly solves an ODE system starting in the point the last event detection yielded. Thus we can only assume that the algorithm produces a path with any desired accuracy if also the event detection locates the hitting times with any desired accuracy. To repeat the purpose of this study, it is precisely this point we address and present for a large class of algorithms a thorough numerical analysis in terms of pathwise convergence and the conditions this presupposes.

Secondly, we note that there exist different but equivalent set-ups for the event detection problem \eqref{sys2}: one as stated in \cite{Zeiseretal} and a second, which is essentially analogous to that, we briefly introduce now. Instead of manipulating equation \eqref{next_jump_time} to obtain \eqref{definition_of_w} we differentiate the survivor function $S$ with respect to $\tau$ and obtain
\begin{equation*}
 \dot{S}(\tau,x)=-S(\tau,x)\, \lambda(\phi(\tau,x))\,.
\end{equation*}
This yields instead of \eqref{sys2} an IVP to solve given by
\begin{equation}\label{sys1}
 \left(\!\!\! \begin{array}{c} \dot{y}\\ \dot{\theta} \\ \dot{S}\end{array}\!\!\!\right)=\left(\!\!\!\begin{array}{c} f(y,\theta)\\ 0 \\ -S\, \lambda(y,\theta)\end{array}\!\!\!\right),\quad \left(\!\!\!\begin{array}{c} y(0)\\ \theta(0)\\ S(0)\end{array}\!\!\!\right)=\left(\!\!\begin{array}{c} X(t_n) \\ 0 \end{array}\!\!\right),\ t\in[0,\infty)\,.
\end{equation}
This system is solved until the component $S(t)$ hits the threshold given by the random variable $U_{2n-1}$. Once again the hitting time is exactly the time until the next jump.
Clearly the two systems \eqref{sys2} and \eqref{sys1} are equivalent in the sense that the $y$ and $\theta$  components of their respective solutions coincide as do the hitting times of the respective thresholds for a given realisation of $U_{2n-1}$.

However, with respect to theoretical analysis the set-up \eqref{sys2} is the most feasible and also we conjecture that for actual implementations it is more efficient as \eqref{sys2} is a simpler type of ODE system than \eqref{sys1}. For this reason, we consider in what follows IVP \eqref{sys2} and denote its solution with respect to the initial value $(x_0,0)$ by $\psi(t,x_0)$. Concerning the different notations introduced so far we summarise
\begin{equation*}
 \psi(t,x_0)=\left(\!\begin{array}{c} \phi(t,x_0) \\ w(t,x_0) \end{array}\!\right) = \left(\!\begin{array}{c} y(t,x_0) \\ \theta_0 \\ w(t,x_0)\end{array}\!\right),\quad \left(\!\!\begin{array}{c} x_0 \\ 0 \end{array}\!\!\right) = \left(\!\!\begin{array}{c} y_0 \\ \theta_0 \\ 0\end{array}\!\!\right)
\end{equation*}
as we make use all of these notations whenever brevity or clarity demands.\medskip

Next, we denote by $\nhat \psi(t,x,h)$, $t\geq 0$, a continuous approximate solution to the IVP \eqref{sys2} with initial condition $(x,0)$ obtained by a continuous numerical ODE method with step size $h$. Consistently, the components of $\nhat\psi$ are denoted by $\nhat\phi$ and $\nhat w$, respectively. An example of a continuous method is the \emph{continuous trapezoidal rule}, which applied to \eqref{ode_set} takes the form
\begin{equation}\label{trapez_1}
\nhat y((n\!+\!1)h)\,=\,\nhat y(nh)+h\tfrac{1}{2}f\bigl(\nhat y(nh),\theta\bigr)+h\tfrac{1}{2}f\bigl(\nhat y((n\!+\!1)h),\theta\bigr)
\end{equation}
with $n=0,\ldots,N\!-\!1,\, h=T/N$ and for $\xi\in[0,1]$
\begin{equation}\label{trapez_2}
\nhat y(nh\!+\!\xi h)\,=\,\nhat y(nh)+h\tfrac{1}{2}\xi(2-\xi)\, f\bigl(\nhat y(nh),\theta\bigr)+h\tfrac{1}{2}\xi^2\, f\bigl(\nhat y((n\!+\!1)h),\theta\bigr)\,.
\end{equation}
That is, \eqref{trapez_1} is the usual trapezoidal rule for ODEs supplemented with an interpolation formula \eqref{trapez_2} to obtain an approximation over the intervals $[nh,(n\!+\!1)h]$ between the discrete grid points.
We remark that although for the numerical methods in Section \ref{section_examples} or the trapezoidal rule just described the parameter $h$ denotes an equidistant step size this need not necessarily be the case.
On the one hand, $h$ may denote the maximal step size used in case of variable step size methods, or, on the other hand, a given error tolerance, cf.~a discussion of deterministic event detection in \cite{Shampine}.
Essentially, $h$ is a defining parameter of the method which, if convergent, converges to the exact solution for \mbox{$h\to 0$\,.} In favour of linguistic simplicity we keep referring to $h$ as the step size in the remainder of the paper.
Further, to keep the presentation and notation simple we restrict ourselves to continuous one-step methods. The implementations in Section \ref{section_examples} are all based on continuous Runge-Kutta methods discussed in \cite{BellenZennaro}. However we expect the subsequent results, i.e., Theorem~\ref{conv_thm} and its Corollary \ref{conv_corol}, to remain valid in the case of continuous multi-step methods.

We refer to \cite{BellenZennaro,hairer} for a general discussion of continuous methods and just briefly collect properties of these methods which we always assume to hold when employed to solve a standard ordinary differential equation of the form
\begin{equation}\label{standard_ode}
\dot y=f(t,y)
\end{equation} 
on the interval $[0,T]$. In order to discuss these properties we employ in this paragraph the following notation: We denote by $y(t,x)$ the exact solution of \eqref{standard_ode} with respect to the initial condition $y(0,x)=x$ and $\nhat y(t,x,h)$ denotes the numerical approximation with respect to the same initial condition and step-size $h$. Then, firstly, the approximate solution obtained by a continuous ODE method satisfies a stability-type estimate of the form
\begin{equation}\label{stab_ODE_method}
\max_{t\in[0,T]}|y(t,x)-\nhat y(t,z,h)|\leq \e^{LT}|x-z|+err(T,h)
\end{equation}
with an error function $err(T,h)$ that satisfies $err(T,h)\to 0$ for $h\to 0$. Moreover, the constant $L>0$ and the function $err$ depend only on the right hand side of \eqref{standard_ode} and especially in the case of the right hand side of \eqref{standard_ode} being Lipschitz continuous, cf.~conditions of Theorem \ref{conv_thm}, $L$ can be chosen to be the Lipschitz constant of $f$ and also $err$ depends in this case only on the Lipschitz constant. Further, $err$ does not depend on the initial condition $x\in E$ and without loss of generality we may assume monotonicity for the error, i.e., $err(T_1,h)\geq err(T_2,h)$ if $T_1>T_2$.
%
%
For the trapezoidal rule \eqref{trapez_1}, \eqref{trapez_2} such a function $err$ is given for small enough $h$ by
\begin{equation*}
 err(T,h)=c\bigl(1+L^{-1}\e^{LT}\bigr)h^2
\end{equation*}
with some appropriate constant $c>0$ if $L$ is a Lipschitz constant to $f$ in \eqref{standard_ode}. 
Particularly, the stability condition \eqref{stab_ODE_method} implies that the method is uniformly convergent on $[0,T]$ in the sense that
\begin{equation}\label{conv_ODE_method}
\lim_{h\to 0}\max_{t\in[0,T]} |y(t,x)-\nhat y(t,x,h)|= 0
\end{equation}
for all initial conditions $x$. We say that the method is convergent of order $p$ if $p$ is the largest integer such that for $h\to 0$
\begin{equation}\label{order_ODE_method}
 \max_{t\in[0,T]} |y(t,y)-\nhat y(t,y,h)|=\landau(h^p),
\end{equation}
in which case $err(T,h)=\landau(h^p)$. Finally, we assume that the numerical approximation satisfies a Lipschitz condition with respect to $t$ on [0,T], i.e.,
\begin{equation}\label{num_method_lipschitz}
|\nhat y(t,x,h)-\nhat  y(s,x,h)|\leq C\,|t-s| \quad \forall\,t,s\in[0,T],
\end{equation}
where the Lipschitz constant $C$ is uniform with respect to the initial condition $x$. In particular, these conditions are satisfied by the one-step methods developed in \cite{BellenZennaro} for the IVP \eqref{sys2} with properties of the right hand side as specified in Theorem \ref{conv_thm}.
\medskip

Thus we obtain an approximation $(\nhat X(t))_{t\in [0,T]}$ for the PDMP $(X(t))_{t\in [0,T]}$ using a continuous ODE method by the following algorithm. Here and subsequently we usually omit to denote the dependence of the numerical approximation on the the step size $h$.\medskip

\begin{algorithm}\label{approx_alg} An algorithm simulating an approximation to a PDMP is given by:
\setlength{\leftmargini}{3.5em}
\begin{enumerate}[align=left]
\setlength{\leftmargini}{3.5em}
  \setlength{\labelwidth}{2.5em}
  \setlength{\labelsep}{1.0em}
\item[Step 1.] Fix the initial time $\nhat t_0=0$ and initial condition $\nhat X(0)=x_0=(y_0,\theta_0)\in E$ and set a jump counter $n=1$. 
\item[Step 2.] Simulate a uniformly distributed random variable $U_{2n-1}$ and solve the IVP \eqref{sys2} with initial condition $(\nhat X(\nhat t_{n-1}),0)$ numerically until
\begin{equation*}
\nhat \tau=\inf\,\{t>0: \nhat w(t)=-\log U\}
\end{equation*}
We take $\nhat t_{n}=\nhat t_{n-1}+\nhat \tau$ as the numerical approximation of the next jump time. Hence, we set
\begin{equation*}
\nhat X(t)=\nhat \phi(t-\nhat t_{n-1},\nhat X(\nhat t_{n-1})) \textnormal{ for } t\in(\nhat t_{n-1},\nhat t_{n}).
\end{equation*}
If $\nhat t_{n}\geq T$, stop at $t=T$.
\item[Step 3.] Otherwise, simulate a post-jump value $\nhat \theta_{n}$ for the piecewise constant component according to the distribution of the jump heights given by $\mu(\nhat \phi(\nhat t_{n}-\nhat t_{n-1},\nhat X(\nhat t_{n-1})),\,\cdot\,)$ via the uniformly distributed random variable $U_{2n}$ and set
\begin{equation*}
 \nhat X(\nhat t_{n})=\left(\!\!\begin{array}{c} \nhat \phi(\nhat t_{n}-\nhat t_{n-1},\nhat X(\nhat t_{n-1})) \\ \nhat \theta_{n}\end{array}\!\!\right).
\end{equation*}
\item[Step 4.] Set $n=n+1$ and start again at Step 2.
\end{enumerate}
\end{algorithm}

Just as hybrid models mix stochastic jump models with continuous deterministic models, hybrid algorithms are essentially constructed by combining a simulation algorithm for the discrete stochastic events, i.e., the SSA, and a numerical method for solving differential equations. In the SSA terminology Algorithm \ref{approx_alg} simulates the stochastic events by Gillespie's direct method, hence such an algorithm may also be called \emph{direct hybrid method}, cf.~\cite{Huisinga}.
The above kind of algorithm, though arising naturally from the model problem, is also related to numerical methods developed for JSDEs. In the case of $\lambda\equiv const.\,$ Algorithm \ref{approx_alg} is exact and turns into jump-adapted methods for JSDEs, see, e.g., \cite{Platen} for a discussion of jump-adapted Taylor methods, in particular for the pure jump case.
We note that derivation or analysis of numerical methods for PDMPs along the lines of It\^o-Taylor expansions employed in the JSDE case is not possible as general PDMPs lack a representation as a solution of a stochastic differential equation. Finally, Algorithm \ref{approx_alg} obviously reduces to the SSA or a purely deterministic ODE method in the degenerate case if either the PDMP is piecewise constant, i.e., $f\equiv const.$, or there are no jumps present at all, i.e., $\lambda\equiv 0$.
\medskip

The following theorem states conditions on the functions $f$ and $\lambda$, which are, together with a valid numerical method for solving \eqref{sys2}, sufficient to guarantee almost sure convergence in the sense of Definition \ref{def_convergence}. The measure $\mu$ is assumed to be as specified in Section \ref{section_PDMP} and  the properties discussed in Section \ref{section_simulations} hold. 

%
%
%
%
%

\begin{theorem}\label{conv_thm} Let $(X(t))_{t\in [0,T]}$ be a regular PDMP with phase space $E$. Let $f$ be bounded, Lipschitz continuous and continuously differentiable on $E$, i.e., there exist constants $M, L$ such that for all $(y,\theta),(z,\vartheta)\in E$ it holds
 \begin{eqnarray}
 |f(y,\theta)|&\,\leq\,& M\,, \label{thm_bound}\\[2ex]
  |f(y,\theta)-f(z,\vartheta)|&\,\leq\,& L\,|(y,\theta)-(z,\vartheta)|\,. \label{thm_lip_cond}
 \end{eqnarray}
Further, assume that the jump rate $\lambda$ is bounded, bounded away from zero, Lipschitz continuous and
continuously differentiable, i.e., there exist constants $\lambda_{\min}, \lambda_{\max}, L$ such that for all $(y,\theta),(z,\vartheta)\in E$ it holds
 \begin{equation}\label{bounded_rate}
 0\,<\,\lambda_{\min}\,\leq\, \lambda(y,\theta)\, \leq\, \lambda_{\max}\,<\,\infty\,,
\end{equation}
\vskip -0.2cm

\begin{equation}\label{thm_rate_lip}
  |\lambda(y,\theta)-\lambda(z,\vartheta)|\ \leq\ L\,|(y,\theta)-(z,\vartheta)|\,.
\end{equation}
Then the algorithms \textnormal{\ref{exact_alg}} and \textnormal{\ref{approx_alg}} are well defined and the numerical approximation $(\nhat X(t))_{t\in[0,T]}$ defined by \textnormal{Algorithm \ref{approx_alg}} converges almost surely to the exact PDMP constructed by \textnormal{Algorithm \ref{exact_alg}} for all initial values \mbox{$X(0)=x_0\in E$} in the sense of \textnormal{Definition \ref{def_convergence}.} Moreover, the continuous components of the PDMP and its approximation converge almost surely uniformly for $h\to 0$, i.e.,
\begin{equation}\label{sup_phase_error_est}
\lim_{h\to 0} \sup_{t\in [0,T]} |Y(t,\omega)-\nhat Y(t,h,\omega)|= 0
\end{equation}
for almost all $\omega\in\Omega$. If, in addition, the continuous ODE method is of order $p$ then also the asymptotic behaviour of \eqref{max_phase_conv}--\eqref{max_time_conv} and \eqref{sup_phase_error_est} is of order $p$.
\end{theorem}

\begin{figure}
\centering
\includegraphics[width=10cm, clip=true, trim=60 75 80 100]{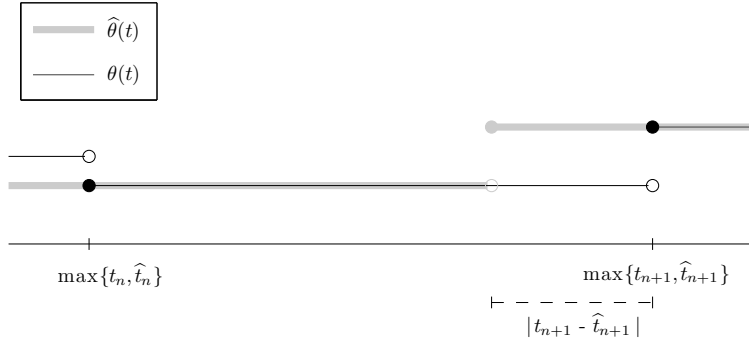}
\caption{Even if the step size $h$ is chosen small enough such that 
the post jump values of the piecewise constant component of the exact PDMP and its approximation coincide, we have that $\sup_{t\in[\max\{t_n,\nhat t_n\},\max\{t_{n+1},\nhat t_{n+1}\}]} |\theta(t)-\nhat \theta(t)| = |\theta(t_n)-\theta(t_{n+1})|>0$ as in general $|t_{n+1}-\nhat t_{n+1}|>0$ for all $h>0$.}\label{disc_traj}
\end{figure}

%
%
%
%
%
%

We briefly comment on the nature of the conditions on $f$ and $\lambda$ in the theorem.
First of all, note that without loss of generality we can choose the same Lipschitz constant in \eqref{thm_lip_cond} and in \eqref{thm_rate_lip}. Particularly this implies that the right hand side of \eqref{sys2} satisfies a Lipschitz condition with constant $L$, hence the global existence of a unique solution is guaranteed. Moreover, this solution is Lipschitz continuous with respect to the initial condition and, as $f,\lambda$ are continuously differentiable, also differentiable with respect to the initial condition \cite{Gruene}.
Secondly, the assumption that $\lambda$ is bounded away from zero is necessary in order for Step 2 in Algorithm \ref{exact_alg} to be well-posed which presupposes, see \cite{Shampine}, on the one hand, that the derivative of the event function is non-zero at the event location, i.e.,
\begin{equation*}
\frac{\dif}{\dif \tau}\int_0^\tau \lambda(\phi(t,x))\,\dif t = \lambda(\phi(\tau,x))\neq 0\,. 
\end{equation*}
On the other hand, it is further necessary for the right hand side of \eqref{sys2} to be bounded in a ball around the event location in phase space, that is, in a ball around the point $z\in E\times\rr_+$ which the solution $\psi$ attains the moment its last component $w$ hits the threshold. 
However, for the stochastic problem we are considering this event can occur at any point in phase space as the initial condition of the IVP \eqref{sys2} as well as the threshold $-\log U_{2n-1}$ are random. Therefore the global boundedness of the right hand side of \eqref{sys2}, i.e., the global boundedness of $f$ and $\lambda$, is required.

We have already mentioned the connection to a convergence result for deterministic event detection \cite{Shampine}. Therein the authors prove that continuous methods conserve their order for problems where one jump event has to be approximated by event detection algorithms over the approximation interval. This result, however, is extended by Theorem \ref{conv_thm} in two ways. First, random instead of deterministic thresholds need to be considered for the simulation of PDMPs and secondly, jump times have to be calculated serially thus a simulation algorithm consists of a sequence of random event detection problems that take the result of the previous event detection as initial condition for the next one. Hence, for an error analysis we have to provide an analysis of the way the single errors of each event detection problem accumulate and are brought forward to a global level.

%
%

\section{Proof of Theorem \ref{conv_thm}}\label{section_proof}
%
%

We begin by recalling some notation which is used throughout the proof.
Let the finite interval $[0,T]$ be fixed. We also fix an $\omega\in\Omega$ such that $N=N(T,\omega)<\infty$, i.e., the realised path $(X(t,\omega))_{t\in[0,T]}$ obtained by Algorithm \ref{exact_alg} has only finitely many jumps in the interval $[0,T]$. By definition of a regular PDMP this includes all $\omega$ except a potential set of measure zero. By $t_n(\omega)$, $n\geq 1$, we denote the jump times of the exact PDMP.
Further, for a fixed step size $h$ of the continuous ODE method used to solve the system \eqref{sys2} we denote by  $(\nhat X(t,h,\omega))_{t\in[0,T]}$ the resulting path of the continuous approximation obtained by Algorithm \ref{approx_alg} with jump times denoted by $\nhat t_n(h,\omega),\ n\geq 1$. 
Finally, we denote by $U_{2n-1}(\omega), U_{2n}(\omega),\ n\geq 1$, the realisations of the i.i.d.~standard uniform random variables defining the $n$th inter-jump time and the $n$th post-jump value, respectively. In the remainder of the proof we generally omit the dependency of any random variable on $\omega$, as well as the dependency of the approximation on $h$.
\medskip

First, we briefly discuss the well-definedness of both Algorithms \ref{exact_alg} and \ref{approx_alg} are well-defined, i.e., we discuss the existence and uniqueness of trajectories for almost all $\omega\in\Omega$, which is inferred immediately from conditions \eqref{thm_lip_cond}--\eqref{thm_rate_lip}. The Lipschitz conditions \eqref{thm_lip_cond} and \eqref{thm_rate_lip} guarantee the existence of a unique global solution $\psi(t,x)$ to the IVP \eqref{sys2} for all initial conditions $x\in E$ and the existence of a numerical approximation $\nhat\psi(t,x)$ follows from the fact that the continuous ODE method is well-defined. Moreover, \eqref{bounded_rate} implies that $w(t,x)$ given by \eqref{definition_of_w} is strictly increasing with $w(0,x)=0$ and $w(t,x)\to\infty$ for $t\to\infty$ and hence \eqref{next_jump_time} has a unique solution for any $U_{2n-1}\in(0,1)$. For the corresponding property of the numerical approximation we note that a method may not produce continuous solutions that are increasing. However, the embedded discrete methods do produce strictly increasing and diverging iterates, which, e.g., immediately follows from the first order condition for Runge-Kutta methods \cite{BellenZennaro}, and condition \eqref{bounded_rate}. 
By the continuity of $\nhat w$ it follows that there exists a finite hitting time to every finite threshold. Uniqueness of this hitting time follows by the use of the generalised inverse in Step 2 of Algorithm \ref{approx_alg}. Moreover, we remark that the same conditions imply that the approximate hitting times $\nhat t_n$ are each  bounded in the step size $h$ for a fixed realisation of the random threshold.

\medskip

We continue by explaining the basic idea of the proof before we present the arguments for the individual steps in detail. As usual for the convergence analysis of numerical methods we compare the approximation and the exact PDMP at discrete time points in $[0,T]$, cf.~the definition of the global error \eqref{max_phase_conv}. Analogously to a standard approach in the convergence analysis we consider a `local' error and show how the error is transported from the local to the global level. Then convergence follows from consistency and numerical stability of the method. However, the `local' error we consider is not the error resulting over one discrete time step in the continuous ODE method but the error that arises over one interval between two successive jumps. This `local' error corresponds to the global error of the ODE method over one inter-jump interval.

Recall that we use a discretisation of the interval $[0,T]$ with grid points given by the jump times $t_n$ and $\nhat t_n$, $n=0,\dots,N$, for the exact and the approximate PDMP, respectively. Therefore the time discretisation for the two processes we have to compare differ in the location of their grid points as well as possibly in the number of grid points, as in general $t_n\not=\nhat t_n$ almost surely and also $T<\nhat t_N$ may be possible with positive probability for some $h>0$.
The former is not a problem as we relate the $n$th inter-jump interval of the exact PDMP to the $n$th inter-jump interval of the approximation regardless of the interval endpoints, see Fig.~\ref{relation}. For this to be well defined the latter, however, is a problem as we thus need the same number of inter-jump intervals in each approximation as there are in the exact solution. Nevertheless, this problem is easily overcome, as we can - without loss of generality - assume that Algorithm \ref{approx_alg} is stopped only after it had at least $N$ jumps. In this way we obtain an approximation with the same number of grid points and which is defined over an interval including $[0,T]$ (at least for small enough $h$, cf.~Section \ref{proof_sec_endpoint_conv}).

\medskip

The proof is organised as follows. In a first step of the proof we derive an error bound for the error in the approximation over one inter-jump interval, i.e., the `local' error of our approximation method. To obtain such a bound we first need a bound on the error in the time until the next jump, i.e., a bound for
\begin{equation}\label{first_def_local_error}
|(t_{n+1}-t_{n})-(\nhat t_{n+1}-\nhat t_{n})|\,.
\end{equation}
Two distinct causes introduce errors in this approximation. On the one hand, the error due to the fact that we numerically approximate the hitting times $t_{n+1}-t_n$ by $\nhat t_{n+1}-\nhat t_n$ and, on the other hand, due to perturbations in the initial conditions as in general $X(t_n)\not=\nhat X(\nhat t_n)$. Therefore, we first consider in Section \ref{proof_first_part} consistency and numerical stability with respect to perturbations in the initial data for the hitting time approximation. Then, in Section \ref{sec_local_phase_err} we consider the `local' error of the method in the PDMP's phase space, whereas in Section \ref{sec_disc_global_err} we prove the first global convergence result, i.e., the limits \eqref{max_phase_conv} and \eqref{max_time_conv}.
The convergence at the interval endpoint, i.e., \eqref{endpoint_error}, is proven in Section \ref{proof_sec_endpoint_conv}.
Finally, in Section \ref{sec_cont_global_err} we extend the discrete convergence result to a continuous result for the PDMP's continuous component, i.e., we prove the convergence \eqref{sup_phase_error_est}.

\subsection{Consistency of the hitting time approximation}\label{proof_first_part}

In this first part of the proof we consider \eqref{first_def_local_error}, i.e., the error in the time until the next jump of the exact PDMP started in $x$ to the time until the next jump obtained by the numerical solution of \eqref{sys2} started from a perturbed initial value $\nhat x$.
Note that we consider the PDMP and its approximation just on one inter-jump interval and hence the components $\theta(t),\,\nhat\theta(t)$ remain constant, however possibly distinct.

To keep the presentation simple we omit the index $n$ in the following, only to be reintroduced in Section \ref{sec_disc_global_err} when we consider the global error of the method. Further, without loss of the generality we set for the exact and approximate inter-jump interval the left endpoint to $0$ and denote by $t$ and $\nhat t$ the jump time of the exact PDMP and its approximation, respectively. Hence the error \eqref{first_def_local_error} reduces to $|t -\nhat t|$ and the two times satisfy the equations
\begin{equation*}
 -\log U = w(t,x)\quad\textnormal{and}\quad -\log U = \nhat w(\nhat t,\nhat x)\,.
\end{equation*}
That is, $t$ and $\nhat t$ are the respective hitting times of the random threshold $-\log U$ of the component $w$ of the exact and approximate solution of the IVP \eqref{sys2}.

\medskip

To estimate $|t-\nhat t|$ we introduce the additional auxiliary time $\ntilde t$ as the time until the next jump of the exact process started in the perturbed initial value $\nhat x$ which is given by
\begin{equation*}
-\log U=w(\ntilde t, \nhat x)\,.
\end{equation*}
Then a simple application of the triangle inequality yields the initial estimate
\begin{equation}\label{consistency_aim}
  |t-\nhat t | \leq |t-\ntilde t|+|\ntilde t-\nhat t|\,.
\end{equation}
Here the first term in the right hand side measures the error in the exact hitting time with respect to perturbations in the initial data whereas the second denotes the numerical error of the method without perturbations in the initial condition. We continue estimating the two terms separately and start with the latter.
We introduce $r(t,h,x)$ as the error of the continuous method's last component defined by 
\begin{equation*}
w(t,x)=\nhat w(t,x)+r(t,h,x)
\end{equation*}
and due to \eqref{stab_ODE_method} it is valid that $|r(s,h,x)|\leq err(t,h)$ for all $s\leq t$ and all $x\in E$.
By definition of the hitting times and the Mean Value Theorem it follows that
\begin{eqnarray*}
0\ \,=\ \, \nhat w(\nhat t,\nhat x)+\log U &=& w(\nhat t,\nhat x)-r(\nhat t,\nhat x,h)+\log U\\[1ex]
&=& w(\ntilde t,\nhat x)+(\nhat t-\ntilde t)\,\nabla_{\!t} w(\vartheta,\nhat x)-r(\nhat t,\nhat x, h)+\log U
\end{eqnarray*}
for some $\vartheta\in[\min\{\ntilde t,\nhat t\},\max\{\ntilde t,\nhat t\}]$. As $w(\ntilde t,\nhat x)+\log U=0$, this equality is equivalent to
\begin{equation*}
(\ntilde t-\nhat t)\,\nabla_{\!t} w(\vartheta,\nhat x)=-r(\nhat t,\nhat x,h),
\end{equation*}
and we obtain due to the definition of $w$ \eqref{definition_of_w} and condition \eqref{bounded_rate} the estimate
\begin{equation}\label{result_one}
|\ntilde t-\nhat t|\leq \frac{|r(\nhat t,\nhat x,h)|}{\lambda_{\min}}\leq \frac{err(\nhat t,h)}{\lambda_{\min}}\,.
\end{equation}

Next we have to consider the term $|t-\ntilde t|$ in the right hand side of \eqref{consistency_aim}. To bound this term we aim for an estimate on the distance in time by the distance of the initial conditions in phase space of pairs $(t,x), (\ntilde t,\nhat x)$ satisfying $w(t,x)=w(\ntilde t,\nhat x)$, i.e., an estimate of the form
\begin{equation}\label{result_two}
 |t-\ntilde t|\leq C_1\,|x-\nhat x|
\end{equation}
for some positive random constant $C_1$ depending only on $t$ and $U$. 
Due to the assumptions in Theorem \ref{conv_thm} the solution of the IVP \eqref{sys2} depends differentiably on the initial value. Thus, proceeding analogously as for error estimate \eqref{result_one} we obtain by the Mean Value Theorem that
\begin{eqnarray*}
0&=&w(\ntilde t,\nhat x)-w(t,x)\\[1ex] &=& w(t,x)+\nabla_{\!t} w(\vartheta,\zeta)\,(\ntilde t-t) + \nabla_{\!x} w(\vartheta,\zeta)\cdot (\nhat x-x)-w(t,x),
\end{eqnarray*}
for some mean values $\vartheta\in[\min\{\ntilde t,\nhat t\},\max\{\ntilde t,\nhat t\}]$, $\zeta\in\{y\,: y=sx+(1-s)\ntilde x,\, s\in[0,1]\}$. Then the Cauchy-Schwarz inequality yields
\begin{equation}\label{result_two_intermediate}
 |t-\ntilde t|\leq \Big|\frac{\nabla_{\!x} w(\vartheta,\zeta)}{\lambda(\phi(\vartheta,\zeta))}\Big|\ |x-\nhat x|.
\end{equation}
This gives an estimate of the form \eqref{result_two} if we bound the pre-factor in the right hand side of \eqref{result_two_intermediate} uniformly.
As the right hand side of the IVP \eqref{sys2} is Lipschitz continuous in $y$ and $\theta$ with constant $L$ it follows that $\phi(t,x)$ depends Lipschitz continuously on the initial condition $x$ with Lipschitz constant $\e^{Lt}$. Moreover, as by the assumptions of the theorem the solution is continuously differentiable with respect to the initial condition it follows that the derivatives satisfy
\begin{equation*}
|\nabla_{\!x} w(s,x)|\leq \e^{Ls} \quad\forall\, x\in E,\ s\geq 0\,. 
\end{equation*}

Hence we obtain for the pre-factor in the right hand side of \eqref{result_two_intermediate} the estimate
\begin{equation*}
\Big|\frac{\nabla_{\!x}w(\vartheta,\zeta)}{\lambda(\phi(\vartheta,\zeta))}\Big|\leq \frac{\e^{\,L\max\{\ntilde t,t\}}}{\lambda_{\min}} \quad\forall\,\vartheta,\,\zeta\,.
\end{equation*}
Finally, it is easy to obtain an a priori bound for $\max\{t,\ntilde t\}$. Obviously, it is valid that
$\max\{\ntilde t,t\}\leq t+|\ntilde t-t|$ and as $t$, $\ntilde t$ are hitting times of the exact solution their difference is globally bounded due to \eqref{bounded_rate}: A straightforward calculation shows that
$|t-\ntilde t|\leq \frac{\lambda_{\max}-\lambda_{\min}}{\lambda_{\max}\lambda_{\min}}\, (-\log U)$ independently of the initial conditions $x$ and $\nhat x$. These estimates applied to \eqref{result_two_intermediate} now yield
\begin{equation}\label{result_two_final}
|t-\ntilde t|\leq \frac{\e^{Lt}}{\lambda_{\min}} U^{-\delta}\,|x-\nhat x|,
\end{equation}
where $\delta=L \frac{\lambda_{\max}-\lambda_{\min}}{\lambda_{\max}\lambda_{\min}}$ and hence we have arrived at an estimate of the form \eqref{result_two}.\medskip

Overall an application of the estimates \eqref{result_one} and \eqref{result_two_final} to the right hand side of \eqref{consistency_aim} yields the stability estimate
\begin{equation}\label{loc_err_time}
 |t - \nhat t| \leq\frac{1}{\lambda_{\min}} \Bigl(\e^{Lt} U^{-\delta}\,|x-\nhat x| + err(\nhat t,h)\Bigr)\,.
\end{equation}

%
%
%
%
%
%

\subsection{The local error in phase space}\label{sec_local_phase_err}

Next we derive a bound of the error over one inter-jump interval in the phase space $E$. We denote the post-jump values of the exact PDMP and its solution, started at initial conditions $x$ and $\nhat x$, after one inter-jump interval by $X(t)$ and $\nhat X(\nhat t)$, respectively. Hence, given a realisation of a standard uniformly distributed random variable $U'$ -- here we use the notation $U'$ to distinguish this random variable from $U$ used defining the inter-jump time -- the `local' error in the phase space is given by
\begin{equation}\label{phase_space_loc_err}
|X(t)-\nhat X(\nhat t)|=|\phi(t,x)+H(\phi(t,x),U')-\nhat\phi(\nhat t,\nhat x)-H(\nhat\phi(\nhat t,\nhat x),U')|\,.
\end{equation}
Note here $H(\,\cdot\,,\,\cdot\,)$ is used as an abbreviation of the correct but cumbersome notation $(0,H(\,\cdot\,,\,\cdot\,))$. To bound this error, we first estimate the difference in continuous components in \eqref{phase_space_loc_err}. Thus, as a first estimate the triangle inequality yields
\begin{equation}\label{cont_local_error}
 |\phi(t,x)-\nhat\phi(\nhat t,\nhat x)|\,\leq\, |\phi(t,x)-\phi(\nhat t,\nhat x)|+|\phi(\nhat t,\nhat x)-\nhat\phi(\nhat t,\nhat x)|.
\end{equation}
The second term in the right hand side of \eqref{cont_local_error} is bounded by $err(\nhat t,h)$ due to stability of the ODE method \eqref{stab_ODE_method}. Another application of the triangle inequality to the first term in the right hand side of \eqref{cont_local_error} yields
\begin{equation}\label{the_rhs}
 |\phi(t,x)-\phi(\nhat t,\nhat x)|\leq |\phi(t,x)-\phi(t,\nhat x)|+|\phi(t,\nhat x)-\phi(\nhat t,\nhat x)|.
\end{equation}
As the solution of the IVP \eqref{sys2} depends Lipschitz continuously on the initial condition the first term in the right hand side of \eqref{the_rhs} satisfies
\begin{equation*}
 |\phi(t,x)-\phi(t,\nhat x)|\leq \e^{L t}\,|x-\nhat x|.
\end{equation*}

To estimate the second term in the right hand side of \eqref{the_rhs} we employ the boundedness \eqref{thm_bound} of $f$ which yields
\begin{equation}
 |\phi(t,x)-\phi(\nhat t,x)|\ \leq \ M\, |t-\nhat t|\,. \label{sec_4_2_analog_est}
\end{equation}

Therefore we overall obtain for the left hand side of \eqref{the_rhs} the estimate
\begin{equation*}
 |\phi(t,x)-\phi(\nhat t,\nhat x)|\leq \e^{Lt}\ |x-\nhat x|+ M\, |t-\nhat t|.
\end{equation*}
Thus employing estimate \eqref{loc_err_time} we get an error bound for the left hand side of \eqref{cont_local_error} by
\begin{equation}\label{local_cont_error}
|\phi(t,x)-\nhat\phi(\nhat t,\nhat x)|\leq \e^{Lt}\bigl(1+\lambda_{\min}^{-1}\, M\,U^{-\delta}\bigr)\,|x-\nhat x|+(1+\lambda_{\min}^{-1}\, M)\,err(\nhat t,h)\,,
\end{equation}
which bounds the error due to the continuous part over one inter-jump interval.
\medskip

Finally, we consider the error introduced by the jump heights $H$ for jumps starting at different points in phase space. In the following we do not bound the error due to the jump heights, but actually show that for sufficiently small step size this error vanishes.
We denote for for all $x,y\in E$ the error by $R(x,y,U')=|H(x,U')-H(y,U')|$ which is finite almost surely and the considerations regarding the structure of the function $H$ in Section \ref{section_simulations} yield that
\begin{equation*}
R(\phi(t,x),\nhat\phi(\nhat t,\nhat x),U')\,=\, 0
\end{equation*}
if $|\phi(t,x)-\nhat\phi(\nhat t,\nhat x)|$ is smaller some bound depending on $\phi(t,x)$ and $U'$. As we have already bound the difference $|\phi(t,x)-\nhat\phi(\nhat t,\nhat x)|$, see \eqref{local_cont_error}, we obtain that $R(\phi(t,x),\nhat\phi(\nhat t,\nhat x),V)=0$ for small enough $h$ and $|\nhat x-x|$.\medskip

To conclude, overall we obtain for the `local' error \eqref{phase_space_loc_err} the estimate
\begin{eqnarray}
|X(t)-\nhat X(\nhat t)|\,\leq\, \e^{Lt}\bigl(1+\lambda_{\min}^{-1}\, M\,U^{-\delta}\bigr)\,|x-\nhat x|+(1+\lambda_{\min}^{-1}\, M)\,err(\nhat t,h)+R(\phi(t,x),\nhat\phi(\nhat t,\nhat x),U')\,,\phantom{xxx}\label{local_error_final}
\end{eqnarray}
where the last term vanishes for small enough $h$ and $|x-\nhat x|$.

%
%
%
%
%
%

\subsection{The discrete global error}\label{sec_disc_global_err}

In order to prove convergence we have to be able to infer from the `local' to the global error. Thus we now reintroduce the indices $n$, hence  $t_n,$, $\nhat t_n$ denote the $n$th jump times of the PDMP and its approximation, respectively, and $X(t_n)$, $\nhat X(\nhat t_n)$ denote the post-jump values of these processes at their jump times. We first show by induction over $n$ that the global error at each jump time, i.e., $|X(t_n)-\nhat X(\nhat t_n)|$ for all $n=1,\ldots,N$, converges to zero for $h\to 0$. Afterward we derive an upper bound for the global error uniform over all jump times to determine the asymptotic order of decrease of the global error.
The splitting of the proof at this point into these two parts is necessary as in order to derive the asymptotic order of convergence we first need to show that the errors $R$ due to different jump heights vanish for small enough step sizes $h$. However, this is a necessary condition of convergence and we can assume it holds for small enough step sizes once the convergence is established.

\medskip

For the induction basis we set $X(0)=\nhat X(0)$, then \eqref{local_error_final} yields the estimate
\begin{equation*}
 |X(t_1)-\nhat X(\nhat t_1)|\ \leq\ 
(1+\lambda_{\min}^{-1} M)\, err(\nhat t_1,h)
+R(\phi(t,x_0),\nhat\phi(\nhat t,x_0),U_2)\,.
\end{equation*}
Obviously for $h$ small enough the last term on the right hand side vanishes due to the convergence of the deterministic method. Then $h\to 0$ implies that $(1+\lambda_{\min}^{-1} M)\, err(\nhat t_1,h)\to 0$ as the approximate hitting time $t_1=t_1(h)$ is bounded in $h$. Therefore we obtain $\lim_{h\to 0}|X(t_1)-\nhat X(\nhat t_1)|=0$\,.\medskip

Next, for the induction step we assume that $\lim_{h\to 0}|X(t_n)-\nhat X(\nhat t_n)|=0$ holds for some $n\geq 1$. Then, for small enough $h$ we have that
\begin{eqnarray*}
|X(t_{n+1})-\nhat X(\nhat t_{n+1})|&\leq& \e^{L(t_{n+1}-t_n)}\bigl(1+\lambda_{\min}^{-1}\,M\,U_{2n+1}^{-\delta}\bigr)\,|X(t_n)-\nhat X(\nhat t_n)|+
(1+\lambda_{\min}^{-1}\,M)\, err(\nhat t_{n+1}-\nhat t_n,h)
\end{eqnarray*}
and hence $\lim_{h\to 0}|X(t_{n+1})-\nhat X(\nhat t_{n+1})|=0$.\medskip

Therefore overall we find that the method is convergent in the sense that
\begin{equation}\label{global_err}
\lim_{h\to 0}\max_{n=1,\dots,N} |X(t_n)-\nhat X(\nhat t_n)| =0\,.
\end{equation}

Analogously we show that also the exact and approximate jump times converge. The error estimate \eqref{loc_err_time} yields
\begin{eqnarray}\label{grid_recursion}
|t_{n+1}-\nhat t_{n+1}|&\,\leq\,& |(t_{n+1}-t_n)-(\nhat t_{n+1}-\nhat t_n)+t_n-\nhat t_n|\nonumber\\[2ex]
 &\leq\,& \e^{L(t_{n+1}-t_n)}\,\lambda_{\min}^{-1}\,U^{-\delta}_{2n+1}\,|X(t_n)-\nhat X(\nhat t_n)|+\lambda_{\min}^{-1}\,err(\nhat t_{n+1}-\nhat t_n,h)+|t_n-\nhat t_n|.\phantom{xxxxx}
\end{eqnarray}
Then \eqref{global_err} and another inductive argument yield that
\begin{equation}\label{global_time_err}
\lim_{h\to 0}\max_{n=1,\dots,N} |t_n-\nhat t_n| = 0\,.
\end{equation}

As these considerations leading to \eqref{global_err} and \eqref{global_time_err} are valid for almost all $\omega\in\Omega$ this completes the proof of the properties  \eqref{max_phase_conv} and \eqref{max_time_conv}. However, to satisfy the condition of pathwise convergence, Definition \ref{def_convergence}, it remains to show the convergence at the interval endpoint \eqref{endpoint_error} which is deferred to Section \ref{proof_sec_endpoint_conv}. At this point we continue to derive estimates for the global errors \eqref{max_phase_conv} and \eqref{max_time_conv} to obtain a rate of convergence.\medskip

To this end we recursively apply \eqref{local_error_final} to its right hand side arriving at
\begin{eqnarray}
|X(t_n)-\nhat X(\nhat t_n)|& \leq& Z_n\,|x_0-\nhat x_0|\\[1ex]
&&\mbox{}+\ Z_n\sum_{k=1}^n Z_k^{-1}\Bigr[(1+\lambda_{\min}^{-1}\, M)\,err(\nhat t_k-\nhat t_{k-1},h)+R(X(t_{k-}),\nhat X(\nhat t_{k-}),U_{2k})\Bigr],\nonumber
\end{eqnarray}
where for $k=1\ldots,N$
\begin{equation}\label{definition_of_Z_n}
 Z_k=\e^{L t_k}\prod_{i=1}^k(1+\lambda_{\min}^{-1}\,M\,U_{2i-1}^{-\delta})\,\geq 1\,.
\end{equation}

As we are considering the asymptotic order of convergence for $h\to 0$ we can assume that
$h$ is small enough such that $R(X(t_{n-}),\nhat X(\nhat t_{n-}),U_{2n})=0$ for all $n=1,\dots,N$. Thus for $x_0=\nhat x_0$ we obtain that for all $n=1,\ldots,N$
\begin{equation}\label{max_phase_error_est}
|X(t_n)-\nhat X(\nhat t_n)|\ \leq\ (1+\lambda_{\min}^{-1}\,M)\, Z_{n}\sum_{k=1}^{n} Z_k^{-1}\,err(\nhat t_k-\nhat t_{k-1},h),
\end{equation}
and therefore it holds that
\begin{equation*}
 \max_{n=1,\dots,N} |X(t_n)-\nhat X(\nhat t_n)|\ \leq\ (1+\lambda_{\min}^{-1}\,M)\,Z_{N}\sum_{k=1}^{N} Z_k^{-1}\,err(\nhat t_k-\nhat t_{k-1},h)\,.
\end{equation*}

Then $err(t,h)=\landau(h^p)$ implies that
\begin{equation}\label{max_phase_error_asym}
\max_{n=1,\dots,N} |X(t_n)-\nhat X(\nhat t_n)|=\landau(h^p)\,.
\end{equation}

Analogously we derive the order of convergence for the jump time approximations. Starting with $t_0=\nhat t_0$ and recursively applying \eqref{grid_recursion} to its right hand side yields 
\begin{equation*}
|t_n-\nhat t_n|
\ \leq\ \lambda_{\min}^{-1}\sum_{k=1}^n\Bigl[err(\nhat t_k-\nhat t_{k-1},h)+\e^{L(t_k-t_{k-1})}\,U_{2k-1}^{-\delta}\ |X(t_{k-1})-\nhat X(\nhat t_{k-1})|\Bigr]\,.
\end{equation*}
Employing \eqref{max_phase_error_est} to estimate the error in phase space yields for small enough $h$ that for all $n=1,\ldots,N$
\begin{equation*}
|t_n-\nhat t_n|\ \leq\ \sum_{k=1}^n \lambda_{\min}^{-1}\Bigl[1+(1+\lambda_{\min}^{-1}\, M)\, Z_{k}^{-1}\sum_{j=k}^{n-1} Z_j\,U_{2j+1}^{-\delta}\,\e^{L(t_{j+1}-t_j)}\Bigr]\,err(\nhat t_k-\nhat t_{k-1},h)\,.
\end{equation*}
Denoting the pre-factor to $err$ in each summand by $W_k^n$ yields that for all $n=1,\ldots,N$
\begin{equation*}
|t_n-\nhat t_n|\leq \sum_{k=1}^{N} W_k^N\,err(\nhat t_k-\nhat t_{k-1},h)
\end{equation*}
and thus $err(t,h)=\landau(h^p)$ implies
\begin{equation}\label{max_time_error_asym}
 \max_{n=1,\dots,N} |t_n-\nhat t_n|=\landau(h^p)\,.
\end{equation}

These considerations leading to \eqref{max_phase_error_asym} and \eqref{max_time_error_asym} are valid for almost all $\omega\in\Omega$ and thus \eqref{max_phase_conv} and \eqref{max_time_conv} follow with the proposed order.

%
%
%
%
%
\subsection{Convergence at the interval endpoint}\label{proof_sec_endpoint_conv}
Finally, for the method to satisfy the definition of pathwise convergence we establish convergence at the interval endpoint $T$. To this end we employ the estimate \eqref{grid_recursion} on the error for the approximation of the $(N+1)$th jump time. This yields
\begin{eqnarray*}
|t_{N+1}-\nhat t_{N+1}|&\leq\,& \e^{L(t_{2N+1}-t_N)}\,\lambda_{\min}^{-1}\,U_{N+1}^{-\delta}\,|X(t_N)-\nhat X(\nhat t_N)|+\lambda_{\min}^{-1}\,err(\nhat t_{N+1}-\nhat t_N,h)+|t_N-\nhat t_N|\,.
\end{eqnarray*}
As all terms in the right hand side of this inequality converge to zero for $h\to 0$ it follows that there exists an $h^\ast>0$, depending on $U_{N+1}$, such that $\nhat t_{N+1}>T$ for all $h<h^\ast$. In this case we then obtain for the error at the interval endpoint the estimate
\begin{eqnarray}
|X(T)-\nhat X(T)|&=&|\phi(T-t_N,X(t_N))-\nhat \phi(T-\nhat t_N,\nhat X(\nhat t_N))|\nonumber\\[1ex]
&\leq& |\phi(T-t_N,X(t_N))-\phi(T-\nhat t_N,X(t_N))|+|\phi(T-\nhat t_N,X(t_N))-\nhat \phi(T-\nhat t_N,\nhat X(\nhat t_N))|\nonumber\\[1ex]
&\leq& M\,|t_N-\nhat t_N|+\e^{L(T-\nhat t_N)}\,|X(t_N)-\nhat X(\nhat t_N)|+err(T-\nhat t_N,h).\label{final_rhs}
\end{eqnarray}
Here we have used in the last inequality on the first term an estimate analogous to \eqref{sec_4_2_analog_est} and on the second term the method's stability estimate \eqref{stab_ODE_method}. Due to previous results and the assumptions on the numerical method the final estimate \eqref{final_rhs} converges to zero for $h\to 0$ for almost all $\omega\in\Omega$ and if the method is of order $p$ then it holds almost surely that
\begin{equation*}
 |X(T)-\nhat X(T)|=\landau(h^p)\,.
\end{equation*}

\subsection{The continuous global error}\label{sec_cont_global_err}

We finish the proof extending the statement \eqref{max_phase_error_asym} to its uniform version \eqref{sup_phase_error_est}. Uniform convergence can only hold for the continuous component $Y(t)$ and cannot hold for the discontinuous component $\theta(t)$, cf.~Fig.~\ref{disc_traj}, which is also shown in a mathematically precise way below.\medskip 

To prove the results in this section we first derive a simple auxiliary result on the asymptotic behaviour of the approximate jump times.
As stated in the last section the convergence of the jump times \eqref{global_time_err} implies that $\nhat t_N>T$ for all $h<h^\ast$ for some $h^\ast$. Moreover there exists an $h^{\ast\ast}<h^\ast$ such that for all $h<h^{\ast\ast}$ it is valid that
\begin{equation}\label{intersection}
 (t_{n-1},t_{n})\cap(\nhat t_{n-1},\nhat t_{n}) \neq \emptyset \quad\forall\ n=1,\dots,N.
\end{equation}
That is, \eqref{intersection} says, that for small enough $h$ the $n$th inter-jump interval of the exact PDMP and its approximation overlap for all $n$. We prove \eqref{intersection} by contradiction. Assume that there exists an $n\leq N$ such that $(t_{n-1},t_{n})\cap(\nhat t_{n-1},\nhat t_{n})=\emptyset$ for all small $h$.  Then it holds that either
\begin{equation*}
 t_{n-1}<t_{n}\leq\nhat t_{n-1} \quad\textnormal{or}\quad \nhat t_{n}\leq t_{n-1}<t_{n}\,,
\end{equation*}
which implies that for all small $h$ either
\begin{equation*}
 0<|t_{n}-t_{n-1}|\leq|\nhat t_{n-1} - t_{n-1}| \quad\textnormal{or}\quad 0<|t_{n-1}-t_{n}|\leq|\nhat t_{n}-t_{n}|.
\end{equation*}
But, both are contradictions to $\lim_{h\to 0}\max_{n=0,\dots,N} |t_n-\nhat t_n|=0$ and thus \eqref{intersection} holds.\medskip

For the remainder of this section we assume that $h<h^{\ast\ast}$, hence $\nhat t_{N+1}>T$ and \eqref{intersection} hold and all errors $R$ vanish. Then, a partition of the time interval $[0,T]$ is given by the discretisation points
\begin{equation*}
 0=:s_0<s_1:=\max\{t_1,\nhat t_1\}<\ldots<s_N:=\max\{t_N,\nhat t_N\}<s_{N+1}:=T\,.
\end{equation*}

It now immediately follows that a uniform convergence result cannot hold for the component $\theta(t)$: On every interval $[s_n,s_{n+1}]$ it is valid that
\begin{equation*}
 \sup_{t\in[s_n,s_{n+1}]} |\theta(t)-\nhat \theta(t)| = |\theta(t_n)-\theta(t_{n+1})|>0 \quad \forall\ h\leq h^{\ast\ast}\,.
\end{equation*}
as $\theta(s_n)=\nhat\theta(s_n)$ and $\theta$ and $\nhat\theta$ jump exactly one time within the interval $[s_n,s_{n+1}]$ to the same post-jump value. By definition, one jumps at $s_{n+1}$ and the other strictly before. Therefore their maximal difference over the interval $[s_n,s_{n+1}]$ is exactly the non-zero jump height independently of $h$.\medskip

After these preliminary considerations we proceed to obtain the continuous convergence result \eqref{sup_phase_error_est}. To this end we use the time grid $s_0,\ldots,s_{N+1}$ and estimate the uniform error of the continuous component $\nhat Y(t)$ on each interval $[s_n,s_{n+1}]$. Note that on each of these intervals it holds that either
\begin{equation*}
\textnormal{i.)}\quad\nhat t_{n+1}< t_{n+1} \quad\textnormal{or}\quad\textnormal{ii.)}\quad t_{n+1}<\nhat t_{n+1}. 
\end{equation*}

We first consider case i.). Recall that we use $y(t,x),\, \nhat y(t,x)$ to denote the first component of the exact and numerical solution of the IVP \eqref{sys2} and we use the symbol $\vee$ to denote the maximum of two expressions. Then
\begin{eqnarray*}
\max_{t\in[s_n,s_{n+1}]}|Y(t)-\nhat Y(t)|&=& \max_{t\in[s_n,\nhat t_{n+1}]}|Y(t)-\nhat Y(t)|\,\vee\,\max_{t\in[\nhat t_{n+1},t_{n+1}]}|Y(t)-\nhat Y(t)|\phantom{xxxxxx}\\
&=&\max_{s\in[0,\nhat t_{n+1}-s_n]}|y(s,X(s_n))-\nhat y(s,\nhat X(s_n))|\\
&& \quad\vee\,  \max_{s\in[0,t_{n+1}-\nhat t_{n+1}]}|y(s,X(\nhat t_{n+1}))-\nhat y(s,\nhat X(\nhat t_{n+1}))|.
\end{eqnarray*}
Both error terms in the right hand side can be estimated using the stability of the continuous ODE method \eqref{stab_ODE_method} which yields
\begin{eqnarray*}
\max_{t\in[s_n,s_{n+1}]}|Y(t)-\nhat Y(t)|&\leq& \Bigl[\e^{L(s_{n+1}-s_n)}\,|X(s_n)-\nhat X(s_n)|\\
&&\ \vee\,\e^{L(t_{n+1}-\nhat t_{n+1})}\,|X(\nhat t_{n+1})-\nhat X(\nhat t_{n+1})|\Bigr]+err(s_{n+1}-s_n,h).
\end{eqnarray*}
As we assume that $h$ is small enough such that $\theta(t_n)=\nhat \theta(\nhat t_n)$ and hence $\theta(s_n)=\nhat\theta(s_n)$ we obtain
\begin{eqnarray*}
\max_{t\in[s_n,s_{n+1}]}\,|Y(t)-\nhat Y(t)|
\!\!\!&\leq&\!\! \Bigl[\e^{L(s_{n+1}-s_n)}|Y(s_n)-\nhat Y(s_n)|\,\vee \\
&&\!\!  \e^{L(t_{n+1}-\nhat t_{n+1})}\bigl(|X(\nhat t_{n+1})-X(t_{n+1})|+|X(t_{n+1})-\nhat X(\nhat t_{n+1})|\bigr)\Bigr]+err(s_{n+1}-s_n,h)\\
&\leq&\!\!  \Bigl[\e^{L(s_{n+1}-s_n)}|Y(s_n)-\nhat Y(s_n)|\,\vee\\
&&\!\!  \e^{L(t_{n+1}-\nhat t_{n+1})}\bigl(M\,|\nhat t_{n+1}-t_{n+1}|+|X(t_{n+1})-\nhat X(\nhat t_{n+1})|\bigr)\Bigr]+err(s_{n+1}-s_n,h).
\end{eqnarray*}
Due to the results \eqref{max_phase_error_asym} and \eqref{max_time_error_asym} it follows that
\begin{equation*}
 \e^{L(t_{n+1}-\nhat t_{n+1})}\bigl(M\,|\nhat t_{n+1}-t_{n+1}|+|X(t_{n+1})-\nhat X(\nhat t_{n+1})|\bigr)+err(s_{n+1}-s_n,h)=\landau(h^p)\,.
\end{equation*}
Hence we obtain the recursive relation
\begin{equation}\label{phase_err_recursion}
 \max_{t\in[s_n,s_{n+1}]}|Y(t)-\nhat Y(t)|\leq \e^{L(s_{n+1}-s_n)}\,|Y(s_n)- \nhat Y(s_n)|+\landau(h^p)\,.
\end{equation}
Secondly, for the case ii.), i.e., $t_{n+1}<\nhat t_{n+1}$, we analogously obtain the estimate
\begin{eqnarray*}
\max_{t\in[s_n,s_{n+1}]}\,|Y(t)-\nhat Y(t)|\!\!\!\!\!&\leq&\!\!\!\! \Bigl[\e^{L(s_{n+1}-s_n)}\,|Y(s_n)-\nhat Y(s_n)|\,\vee\\
&&\!\!\!\! \e^{L(\nhat t_{n+1}- t_{n+1})}\bigl(|X(t_{n+1})-\nhat X(\nhat t_{n+1})|+|\nhat X(\nhat t_{n+1})-\nhat X(t_{n+1})|\bigr)\Bigr]+err(s_{n+1}-s_n,h)\,.
\end{eqnarray*}

To estimate the term $|\nhat X(\nhat t_{n+1})-\nhat X(t_{n+1})|$ we employ the Lipschitz condition \eqref{num_method_lipschitz}. Hence we obtain again a recursive relation of the form \eqref{phase_err_recursion} . Thus such a relation holds for the interval $[s_n,s_{n+1}]$ in either cases.\medskip

Finally, recursively applying \eqref{phase_err_recursion} to its right hand side yields the estimate
\begin{equation*}
 \max_{n=0,\dots,N}\max_{t\in[s_n,s_{n+1}]}|Y(t)-\nhat Y(t)|\leq \e^{LT}|Y(0)-\nhat Y(0)|+\landau(h^p)\,.
\end{equation*}
These considerations hold true for almost all $\omega\in\Omega$, hence the uniform convergence \eqref{sup_phase_error_est} follows for $Y(0)=\nhat Y(0)$ and the proof of Theorem \ref{conv_thm} is completed.

\section{Extensions of the convergence theorem}\label{section_extension}

The results of Theorem \ref{conv_thm} are easily generalised to arbitrary PDMPs with right continuous paths. That is, we now consider PDMPs that do not have any qualitative differences in their components: All components allow for discontinuities, i.e., $\mu$ is a Markov kernel onto $\rr^{d+m}$ -- however still discretely supported --, and in between jumps the trajectories follow a deterministic motion given by an ODE
\begin{equation*}
 \left(\!\!\! \begin{array}{c} \dot{y}\\ \dot{\theta}\end{array}\!\!\!\right)=\left(\!\!\!\begin{array}{c} f_1(y,\theta)\\ f_2(y,\theta)\end{array}\!\!\!\right)\,,
\end{equation*}
i.e., the component $\theta$ is not piecewise constant anymore. Clearly, a continuous convergence result such as \eqref{sup_phase_error_est} cannot hold anymore. However, it is straightforward to see that the result in Theorem \ref{conv_thm} on the asymptotic behaviour of the errors at the jump times, at the interval end point and the errors in the jump times is still valid if $f_1$ and $f_2$ satisfy the conditions specified therein for $f$.\medskip

In a further extension we consider the just discussed general right-continuous structure of a PDMP but assume that the distributions $\mu(x,\cdot)$ of the jump heights  are continuous on $\rr^{d+m}$ for all $x\in E$ instead of discrete. Corollary 23.4 in \cite{Davis2} guarantees as in Section \ref{section_simulations} the existence of a function $H: E\times(0,1)\to \rr^{d+m}$ providing realisations of the random jump heights from realisations of standard uniformly distributed random variables. In order to derive a convergence result in this case we impose on $H$ the following Lipschitz condition: there exists a deterministic constant $L_H<\infty$ such that for all $x,y\in \rr^d$ it holds almost surely that
\begin{equation}\label{jump_errors_disappear3}
 |H(x,U_{2n})-H(y,U_{2n})| \ \leq\ L_H\,|x-y|\quad\forall\,n\geq 1\,.
\end{equation}
Then by only a small change in the proof of Theorem \ref{conv_thm} we are able to establish the following convergence result.

\medskip

\begin{corollary}\label{conv_corol} Let $(X(t))_{t\in[0,T]}$ and $(\nhat X(t))_{t\in[0,T]}$ be a regular PDMP and its approximation, respectively, which have right-continuous paths. The distribution of the jump heights is continuous and \eqref{jump_errors_disappear3} is satisfied. Then, under the conditions of Theorem \ref{conv_thm} Algorithm \ref{approx_alg} converges pathwise in the sense of Definition \ref{def_convergence}, i.e., the global error satisfies \eqref{max_phase_conv}--\eqref{max_time_conv}. Moreover, if the embedded continuous ODE method is of order $p$, then the order of convergence is $p$.
\end{corollary}

%
%
%
%

\medskip

\begin{proof} The proof of \eqref{max_phase_conv}--\eqref{max_time_conv} in case of this general class of PDMPs works analogously to the proof of Theorem \ref{conv_thm} in Section \ref{section_proof}. The arguments employed become even less technical. In particular, the last paragraph in Section \ref{sec_local_phase_err} dealing with the error in the jump heights becomes redundant as does the induction argument in Section \ref{sec_disc_global_err}. We restrict the presentation to point out the main difference.

The proof proceeds as the proof in Section \ref{section_proof} with the main difference in estimating the local error in phase space, cf.~Section \ref{sec_local_phase_err}. Using the Lipschitz condition \eqref{jump_errors_disappear3} we first obtain from \eqref{phase_space_loc_err} the estimate
\begin{equation*}
|X(t)-\nhat X(\nhat t)|=(1+L_H)\,|\phi(t,x)-\nhat\phi(\nhat t,\nhat x)|\,.
\end{equation*}
Then, estimating the difference $|\phi(t,x)-\nhat\phi(\nhat t,\nhat x)|$ as in Section \ref{sec_local_phase_err} yields for the local error in phase space an estimate by
\begin{equation*}
|X(t)-\nhat X(\nhat t)|\, \leq\, \e^{Lt}(1+L_H)\bigl(1+\lambda_{\min}^{-1}\, M\,U^{-\delta}\bigr)\,|\nhat x-x|+(1+L_H)(1+\lambda_{\min}^{-1}\, M)\,err(\nhat t,h)\,,
\end{equation*}
cf.~the local error estimate \eqref{local_error_final}. Starting with $X(0)=\nhat X(0)$ and a recursive application of this last inequality, see Section \ref{sec_disc_global_err}, yields the global error bound 
\begin{equation*}
 \max_{n=1,\dots,N} |X(t_n)-\nhat X(\nhat t_n)|\ \leq\ (1+L_H)(1+\lambda_{\min}^{-1}\,M)\,Z_{N}\sum_{k=1}^{N} Z_k^{-1}\,err(\nhat t_k-\nhat t_{k-1},h)\,.
\end{equation*}
As this calculations are valid for almost all $\omega\in\Omega$ the convergence \eqref{max_phase_conv} follows and is of order $p$ in case $err(t,h)=\landau(h^p)$.
The proofs for the limits \eqref{max_time_conv} and \eqref{endpoint_error} work completely analogous as for Theorem \ref{conv_thm}, see Sections \ref{sec_disc_global_err}  and \ref{proof_sec_endpoint_conv}, respectively.
\end{proof}

\section{Numerical examples}\label{section_examples}

\begin{figure}
\begin{equation*}
\qquad\quad\qquad\begin{array}{cccccc}
\multicolumn{2}{c}{\begin{array}{c|c}
 c & A \\
\hline\\[-2ex]
 & \beta^T
\end{array}} & \multicolumn{2}{c}{\begin{array}{c|c}
 0 & 0 \\
\hline\\[-2.5ex]
 & 1
\end{array}} & \multicolumn{2}{c}{\qquad\quad\begin{array}{c|cc}
 0 & 0 & 0\\
 1 & \frac{1}{2} & \frac{1}{2}\\[0.5ex]
\hline\\[-2ex]
 & \frac{1}{2} & \frac{1}{2}
\end{array}\qquad\quad}\\[5ex]
\multicolumn{2}{l}{} & \multicolumn{2}{l}{\qquad\quad b_1(\xi)=\xi\qquad\quad} & \multicolumn{2}{l}{\qquad\qquad b_1(\xi)=\frac{1}{2}\xi(2-\xi)\qquad\quad} \\[1.5ex]
\multicolumn{2}{l}{} & \multicolumn{2}{l}{} & \multicolumn{2}{l}{\qquad\qquad b_2(\xi)=\frac{1}{2}\xi^2} \\[2ex]
\multicolumn{2}{c}{\textnormal{(a)}} & \multicolumn{2}{c}{\textnormal{(b)}} & \multicolumn{2}{c}{\textnormal{(c)}} \\[3ex]
\multicolumn{3}{c}{\begin{array}{c|cc}
 \frac{1}{3} & \frac{5}{12} & -\frac{1}{12}\\[0.75ex]
 1 & \frac{3}{4} & \frac{1}{4}\\[0.5ex]
\hline\\[-2ex]
 & \frac{3}{4} & \frac{1}{4}
\end{array}} & \multicolumn{3}{c}{\begin{array}{c|ccc}
 0 & 0 & 0 & 0\\
 \frac{1}{2} & \frac{5}{24} & \frac{1}{3} & -\frac{1}{24}\\[0.75ex]
 1 & \frac{1}{6} & \frac{2}{3} & \frac{1}{6}\\[0.75ex]
\hline\\[-2ex]
 & \frac{1}{6} & \frac{2}{3} & \frac{1}{6}
\end{array}}\\[7ex]
\multicolumn{3}{l}{\quad b_1(\xi)=\frac{3}{4}\xi(2-\xi)\quad} & \multicolumn{3}{l}{\qquad\qquad\quad b_1(\xi)=2\xi(\frac{1}{3}\xi^2-\frac{3}{4}\xi+\frac{1}{2})\qquad\qquad\quad}\\[1.5ex]
\multicolumn{3}{l}{\quad b_2(\xi)=\frac{3}{4}\xi(\xi-\frac{2}{3})} & \multicolumn{3}{l}{\qquad\qquad\quad b_2(\xi)=4\xi^2(\frac{1}{2}-\frac{1}{3}\xi)}\\[1.5ex]
\multicolumn{3}{l}{} & \multicolumn{3}{l}{\qquad\qquad\quad b_3(\xi)=2\xi^2(\frac{1}{3}\xi-\frac{1}{4})}\\[2ex]
\multicolumn{3}{c}{\textnormal{(d)}} & \multicolumn{3}{c}{\textnormal{(e)}}
\end{array}
\end{equation*}
\caption{Coefficients and interpolation formulae of the continuous Runge-Kutta methods implemented for the numerical examples: \textnormal{(a)} general Butcher tableau, \textnormal{(b)} forward Euler method, \textnormal{(c)} trapezoidal rule, \textnormal{(d)} RadauIIa method, \textnormal{(e)} LobattoIIIa method.}\label{butcher_tableaus} 
\end{figure}

To illustrate the theoretical findings on the order of convergence we have implemented simulation methods for PDMPs based on continuous Runge-Kutta methods of different order. For an IVP
\begin{equation*}
 \dot y=f(t,y),\quad y(0)=y_0
\end{equation*}
with $y(t),\,y_0 \in\rr^d$ and $t\in[0,T]$ an \emph{$s$-stage continuous Runge-Kutta method} is a discretisation scheme of the form
\begin{equation}
\nhat y((n+1)h)\, =\, \nhat y(nh)+h\sum_{i=1}^s \beta_i f(nh+c_ih,k_i),\quad n=0,\ldots,N-1,\, h=T/N
\end{equation}
with stage values $k_i,\,i=1\ldots,s$, given by
\begin{equation}\label{stage_values}
 k_i=\nhat y(nh)+h\sum_{j=1}^s a_{ij} f(nh+c_jh,k_j)
\end{equation}
and combined with an interpolation formula
\begin{equation}
\nhat y(nh+\xi h)=\nhat y(nh)+h\sum_{i=1}^s b_i(\xi)\,f(nh+c_ih,k_i),\quad 0\leq\xi\leq 1
\end{equation}
for the approximation on the intervals between the discretisation points $nh$.
The coefficients $\beta=(\beta_1,\ldots,\beta_s)^T$, $A=(a_{ij})_{i,j=1,\ldots,s}$ and $c=(c_1,\ldots,c_s)^T$ are given by Butcher tableaus as in Fig.~\ref{butcher_tableaus} (a). We implemented the \emph{explicit Euler method} (order 1), the \emph{trapezoidal rule} (order 2), the \emph{$2$-stage RadauIIa method} (order 3) and the \emph{$3$-stage LobattoIIIA method} (order 4). The coefficients and interpolation polynomials $b_i(\xi)$ for these methods are given in Fig.~\ref{butcher_tableaus} (b)--(e) taken from \cite{BellenZennaro}. Each experiment consists of the same trajectory simulated using Algorithm \ref{approx_alg} based on the different continuous Runge-Kutta methods with decreasing step sizes $h$. We compared the approximations to a reference solution as an exact solution is not available. The reference solution is an approximation of the same trajectory simulated with very high accuracy. We stopped decreasing the step size $h$ for an implementation of Algorithm \ref{approx_alg} when the approximations entered the error range of the reference solution. 
\medskip

%
%

We have applied the methods to a stochastic hybrid version of the Hodgkin-Huxley model for a patch of neuronal membrane \cite{BuckwarRiedler,Wainrib1}. We consider the model for two different parameters sets denoted by (P1) and (P2). The numerical values of both parameter sets can be found in the appendix. The set (P1) is the original Hodgkin-Huxley model for the squid giant axon. Since its introduction in \cite{HodgkinHuxley} this model is the standard example in neuroscience for models of excitable membranes. The second set (P2) is taken from \cite{Minoetal} wherein the authors experimentally compare the performance of various pseudo-exact algorithms with respect to certain test statistics.

The model describes the \emph{electrical potential difference} across a neuronal membrane which is dynamically changed due to varying ionic currents through the membrane. These currents are mediated by membrane proteins that span across the membrane and in changing their shape these open and close pores. Such proteins are called \emph{ion channels} and are named after the class of ions they allow selectively to pass. Changes in channel states occur in a random, memoryless fashion with instantaneous rates depending on the current potential difference. Conversely, the current through the membrane is proportional to the number of open states. For a more physiological explanation and derivation of the mathematical model we refer to \cite{Koch} and \cite{BuckwarRiedler} for its description using PDMPs.

The hybrid version of the Hodgkin-Huxley neuron model is a 14--dimensional PDMP where a one-dimensional continuous variable $Y(t)\in\rr$ models the difference in the electrical potential. The remaining, piecewise constant variables $\theta(t)\in\rr^{13}$ record the states of the ion channels immersed in the membrane. In the Hodgkin-Huxley model there are two different families of ion channels: sodium ($Na^+$) and potassium ($K^+$) channels with $8$ and $5$ distinct states, respectively. Their kinetic schemes are given in Fig.~\ref{channel_kinetics} in the appendix. Each component of $\theta(t)$ thus counts the number of channels in the specific state, i.e., $\theta_1,\ldots,\theta_8$ correspond to the states of the $Na^+$--channels and $\theta_9,\ldots,\theta_{13}$ to the states of the $K^+$--channels.

The characteristics of the PDMP are as follows. Firstly, the family of ODEs \eqref{ode_set} defining the inter jump evolution of the PDMP is given by the equation
\begin{equation}\label{membrane_balance_eq}
C\dot y= -\nbar g_\textnormal{Na}\theta_8(y-E_\textnormal{Na})-\nbar g_\textnormal{K}\theta_{13}(y-E_\textnormal{K})-\nbar g_\textnormal{L}(y-E_\textnormal{L})+I(t)\,.
\end{equation}

Secondly, the instantaneous jump rate $\lambda$ is given by
\begin{eqnarray*}
\lambda\bigl((y,\theta)\bigr)&=&
\left(\!\!\begin{array}{c}
a_m(y) \\ b_m(y) \\ a_h(y) \\ b_h(y)
      \end{array}\!\!\right)^T\,
\left(\begin{array}{cccccccc}
 3 & 2 & 1 & 0 & 3 & 2 & 1 & 0\\
 0 & 1 & 2 & 3 & 0 & 1 & 2 & 3\\
 1 & 1 & 1 & 1 & 0 & 0 & 0 & 0\\
 0 & 0 & 0 & 0 & 1 & 1 & 1 & 1\\
\end{array}\right)\,\left(\!\!\begin{array}{c}\theta_1\\ \vdots \\ \theta_8\end{array}\!\!\right)\\
& &\qquad\qquad\qquad\qquad\qquad\qquad+
\left(\!\!\begin{array}{c}
a_n(y) \\ b_n(y)
      \end{array}\!\!\right)^T\,
\left(\begin{array}{cccccccc}
 4 & 3 & 2 & 1 & 0 \\
 0 & 1 & 2 & 3 & 4 
\end{array}\right)\,\left(\!\!\begin{array}{c}\theta_9\\ \vdots \\ \theta_{13}\end{array}\!\!\right)\,.
\end{eqnarray*}

\begin{figure}
\centering
\includegraphics[width=0.95\textwidth, clip=true, trim=5 0 15 0]{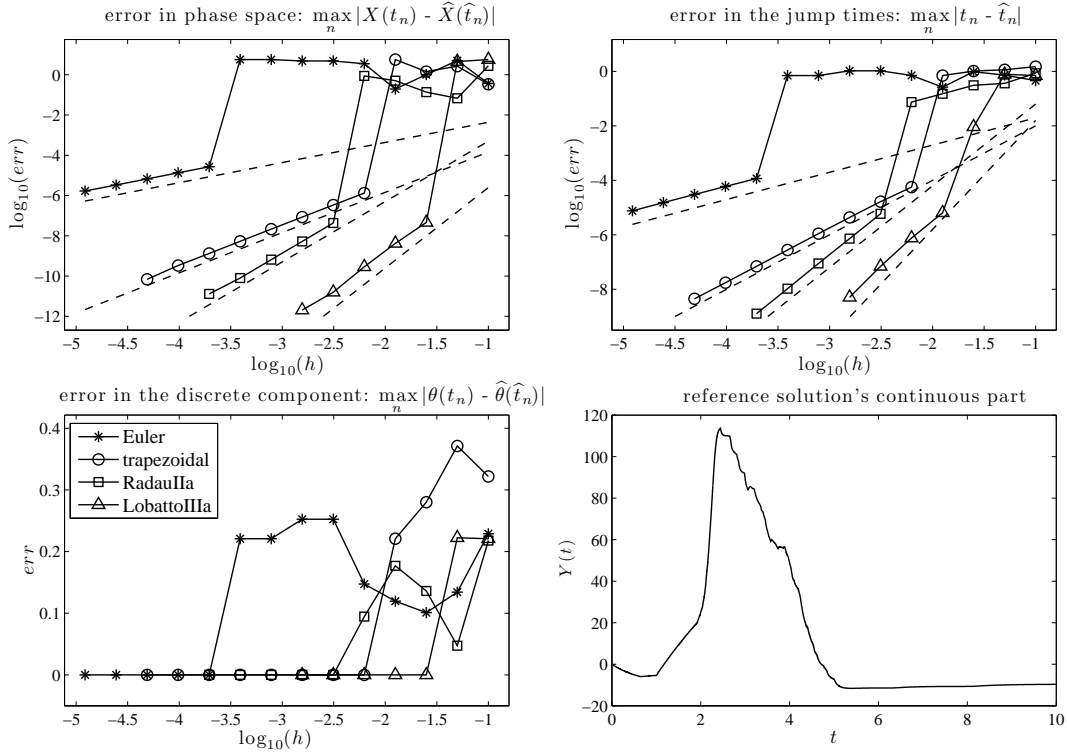}
\caption{A first exemplary numerical experiment for a trajectory to the parameter set \textnormal{(P1)}. We have plotted the errors in phase space (top left), jump times (top right) and in the discrete component (bottom left). The path of the reference solution's continuous component calculated using \textnormal{MATLAB}$\textsuperscript{\textregistered}$'s \textnormal{\texttt{ode45}} is shown in the bottom right panel.}\label{ex_1_traj}
\end{figure}

\begin{figure}
\centering
\includegraphics[width=0.95\textwidth, clip=true, trim=5 0 15 0]{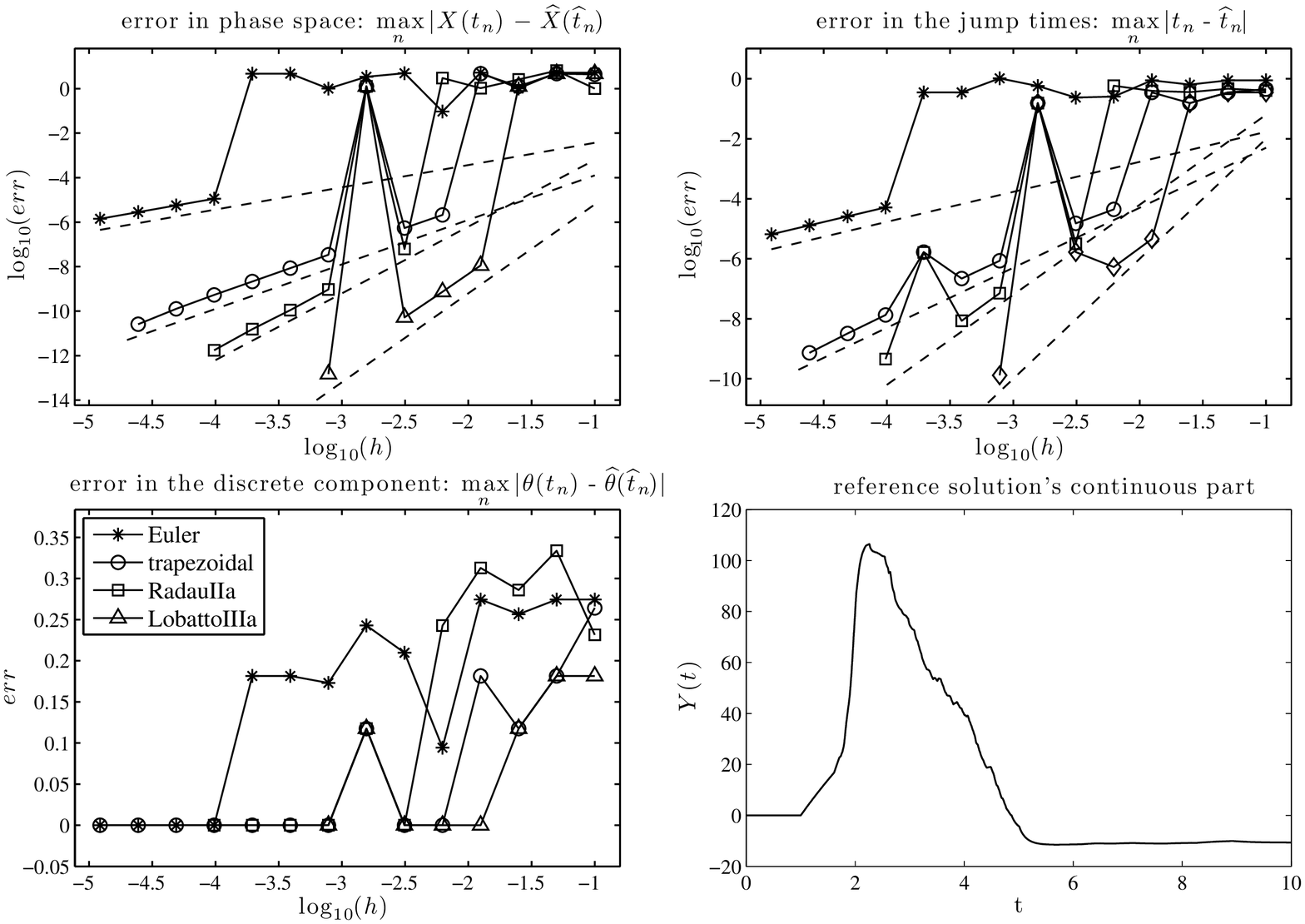}
\caption{A second exemplary numerical experiment for a trajectory to the parameter set \textnormal{(P1)}. We have plotted the errors in phase space (top left), jump times (top right) and in the discrete component (bottom left). The path of the reference solution's continuous component calculated using the presented PDMP method based on the continuous \textnormal{LobattoIIIA} is shown in the bottom right panel.}\label{ex_2_traj}
\end{figure}

Thirdly, we specify the transition measure $\mu$. The point probability of the event that one channel changes from state $i$ to state $j$ is given by the transition rate of state $i$ to state $j$ times the number of channels in state $i$ divided by the total instantaneous rate $\lambda$. For example, the probability of the event of one channel changing from state $1$ to $2$ -- conditional on the jump time being $t$ -- is given by
\begin{equation*} 
\frac{3\,a_m(Y(t))\,\theta_1(t-)}{\lambda\bigl((Y(t),\theta(t-))\bigr)}=: \mu\bigl((Y(t),\theta(t-)),\{(-1,1,0,\ldots,0)^T\}\bigr)\,.
\end{equation*}
All other events have zero probability, that is two channels do not change states simultaneously almost surely.

We remark that \eqref{membrane_balance_eq} differs from the general form \eqref{ode_set} due to the added time dependent function $I(t)$ which denotes the external current input to the system. However, for the numerical experiments we present the input is piecewise constant of the form $I(t)=const.\cdot\mathbb{I}_{(T_1,T_2]}(t)$ (\emph{monophasic input}) and hence one may think of the model as PDMPs `glued' together at the times $T_1,\,T_2$ with the final state being the initial condition of the next. The reason for this initialising time is that in numerical experiments one wants to sample an initial condition from the membrane at rest, i.e., with $I(t)\equiv 0$, for the start of the response to an input at time $T_1$.

\medskip

We present two sample trajectories for the parameter set (P1). The reference solution for the trajectory in Fig.~\ref{ex_1_traj} was calculated using MATLAB$\textsuperscript{\textregistered}$'s \texttt{ode45} implementation with \texttt{AbsTol} / \texttt{RelTol} $=2.22045\e^{-14}$. We used the built-in event detection for the calculation of the hitting times. For the trajectory in Fig.~\ref{ex_2_traj} the reference solution was calculated using the LobattoIIIa method with step size $h=5\e^{-6}$. On the one hand the use of LobattoIIIA illustrates that high order methods yield good results in the absolute error range for reasonable equidistant step sizes.
On the other hand the use of \texttt{ode45}, which employs automated step-size selection, illustrates that methods employing automated step size selection yield good accuracy results and thus can in principle be used.
However, we found that for the same level of accuracy simulation times for the equidistant LobattoIIIA method were considerably shorter than for \texttt{ode45}. Thus an equidistant implementation performs better as a simulation using MATLAB$\textsuperscript{\textregistered}$'s time stepping algorithm contrary to the latter being specifically developed to speed up simulations. We note that these findings are not artefacts of the two specific paths presented but are persistent throughout all the simulations we conducted.

The error plots in the upper panels of Fig.~\ref{ex_1_traj}, i.e., the error in phase space \eqref{max_phase_conv} (upper left panel) and the error in the jump times \eqref{max_time_conv} (upper right panel) illustrate very clearly the theoretical result in Theorem \ref{conv_thm}: For step size $h$ small enough such that the global error in the discrete component vanishes (lower left panel) the theoretical asymptotic order of convergence for the different methods is observed. For a guide of the eye we added
lines of slope $1$, $2$, $3$ and $4$ beneath the errors of the methods with the respective orders.
The plots of the second trajectory shown in Fig.~\ref{ex_2_traj} are a slightly less clear illustration of the asymptotic order of convergence. Here, the discrete error of the higher order methods has already vanished however reappears for a certain step size $h$. However, if we remove this outlier we observe the asymptotic order underlying the approximation as seen by comparison with the straight lines.


As mentioned discussing the reference solution a naive use of automated step-size control with the intention of speeding up simulations and controlling the error can be misleading.
However, it would be a task for an efficient and practicable pathwise step size detection and error control to detect such `bad' step sizes as occur in the second trajectory and either avoid them or minimise their effect. 
These are crucial points for the implementation and performance of the algorithms. They demand for a further thorough investigation in this direction, e.g., an analysis of the shortcomings of the step-size detection as employed in \texttt{ode45}, which we have not attempted to do in this study.

%
%

\begin{figure}
\centering
\includegraphics[width=0.95\textwidth, clip=true, trim=5 0 15 0]{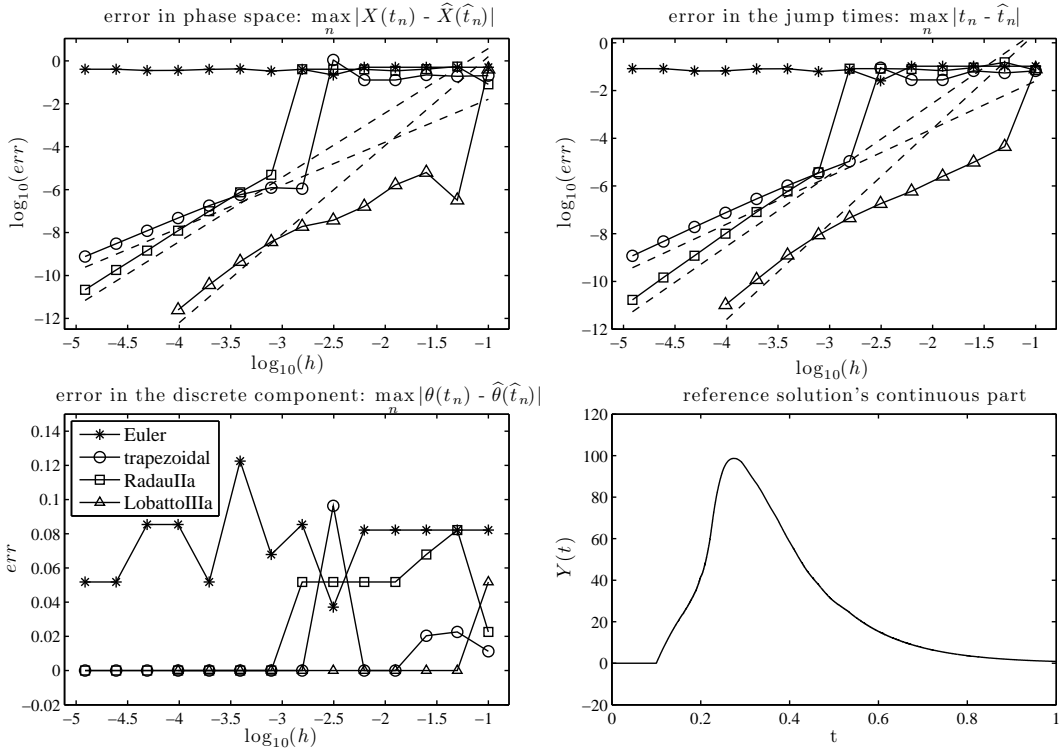}
\caption{An exemplary numerical experiment for a trajectory to the parameter set \textnormal{(P2)}. We have plotted the errors in phase space (top left), jump times (top right) and in the discrete component (bottom left). The path of the reference solution's continuous component calculated using the presented PDMP method based on the continuous \textnormal{LobattoIIIA} is shown in the bottom right panel.}\label{ex_3_traj}
\end{figure}

Finally, in a second example we consider a variant of the Hodgkin-Huxley model
only dealing with currents due to $Na^+$--channels, that is the model reduces to $9$ dimensions as $\theta_9(t),\ldots,\theta_{13}(t)\equiv 0$. The channel density is higher and the current per channel smaller as for the first parameter set which overall renders the trajectories `less noisy', cf.~the sample trajectory of the reference solution in Figs.~\ref{ex_1_traj} and \ref{ex_2_traj} with the one in Fig.~\ref{ex_3_traj}. The reference solution was again calculated using the LobattoIIIA method with step size $h=5\e^{-6}$. Overall we find the same behaviour as in the first example.
For small enough step sizes such that the errors in the discontinuous components vanish, the predicted order of the asymptotic error behaviour is observed. 
However, note that even for the smallest step size considered the Euler method lacks the accuracy to approximate the trajectory correctly. Its error is orders of magnitude larger than the error of higher order Runge-Kutta method with rather large step size. This observation strongly supports the use of higher order methods in practice.

\section*{Acknowledgements:} The author acknowledges support of the EPSRC grant EP/E03635X/1 ``Centre for Analysis and Nonlinear Partial Differential Equations'' and from the BC/DAAD ARC project ``Robust simulation of networks with random switching'' (1349/50021880). The work in this study was accomplished while the author was at Heriot-Watt University, Edinburgh, UK. Further, the author thanks Prof.~E.~Buckwar, her thorough reading and insightful comments helped to improve the presentation in the present paper substantially.

\appendix

\begin{figure}
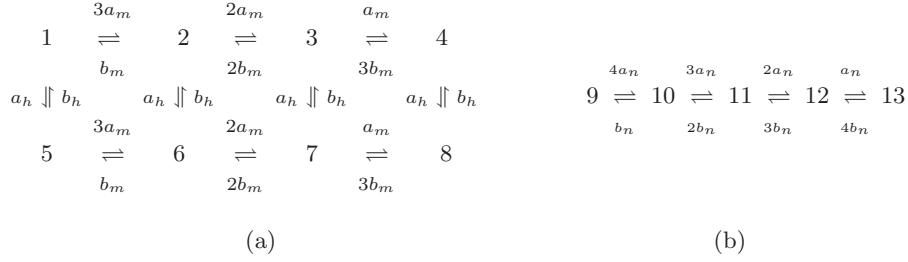

\small{
\begin{equation*}
\left.\begin{array}{cc}

\left.\begin{array}{ccccccc}
      &\!\!\textnormal{\scriptsize $3a_m$}\!\!                 &       &\!\! \textnormal{\scriptsize $2a_m$}\!\!                 &       & \!\!\textnormal{\scriptsize $a_m$}\!\!	                &      \\
\!\!1\!\! &\!\! \rightleftharpoons & 2\!\! &\!\! \rightleftharpoons\!\! &\!\! 3\!\! &\!\! \rightleftharpoons\!\! &\!\! 4\\
      & \!\!\textnormal{\scriptsize $b_m$}\!\!                 &       &\!\! \textnormal{\scriptsize $2b_m$}\!\!                 &       & \!\!\textnormal{\scriptsize $3b_m$}\!\!                 &      \\
\!\!\textnormal{\scriptsize $a_h$}\downharpoonleft\hspace{-0.8mm}\upharpoonright \textnormal{\scriptsize $b_h$}\!\!& &\!\!\textnormal{\scriptsize $a_h$}\downharpoonleft\hspace{-0.8mm}\upharpoonright \textnormal{\scriptsize $b_h$}\!\!& &\!\! \textnormal{\scriptsize $a_h$}\downharpoonleft\hspace{-0.8mm}\upharpoonright \textnormal{\scriptsize $b_h$}\!\! &&\!\! \textnormal{\scriptsize $a_h$}\downharpoonleft\hspace{-0.8mm}\upharpoonright \textnormal{\scriptsize $b_h$}\\
      & \!\!\textnormal{\scriptsize $3a_m$}\!\!                 &       & \!\!\textnormal{\scriptsize $2a_m$}\!\!                 &       & \!\!\textnormal{\scriptsize $a_m$}\!\!                 &      \\
\!\!5\!\! &\!\! \rightleftharpoons\!\! &\!\! 6\!\! &\!\! \rightleftharpoons\!\! &\!\! 7\!\! &\!\! \rightleftharpoons\!\!& \!\!8\!\!\\
      &\!\! \textnormal{\scriptsize $b_m$}\!\!                 &       &\!\! \textnormal{\scriptsize $2b_m$} \!\!                &       & \!\!\textnormal{\scriptsize $3b_m$}\!\!                 &\end{array}\right.

\ \ \
&
\ \ \
\left.\begin{array}{ccccccccc}
&\!\! \textnormal{\tiny $4a_n$}\!\!& &\!\!\textnormal{\tiny $3a_n$}\!\!&&\!\!\textnormal{\tiny $2a_n$}\!\!& & \!\!\textnormal{\tiny $a_n$}\\
9\!\! &\!\!\rightleftharpoons\!\!&\!\! 10\!\! & \!\!\rightleftharpoons\!\! &\!\! 11 \!\!&\!\! \rightleftharpoons\!\! &\!\!12\!\!&\!\!\rightleftharpoons\!\! &\!\! 13\!\!\\
&\!\!\textnormal{\tiny $b_n$}\!\!& &\!\! \textnormal{\tiny $2b_n$}\!\! &&\!\! \textnormal{\tiny $3b_n$}\!\! && \!\!\textnormal{\tiny $4b_n$}\!\!\end{array}\right.

\\ & \\ \textnormal{(a)}  & \textnormal{(b)}
\end{array}\right.
\end{equation*}
\caption{Kinetic scheme of a \textnormal{(a)} $Na^+$ and \textnormal{(b)} $K^+$ ion channel in the standard Hodgkin-Huxley model. The states $8$ and $13$ are the conducting states of channels, respectively. The rates $a_x=a_x(y),\ b_x=b_x(y)$ are dependent on the transmembrane potential $y$.}
}\label{channel_kinetics}
\end{figure}

\section*{Appendix.}\label{coeff_rateFunc}

The following are the rate functions and coefficients for the original Hodgkin-Huxley model of the squid giant axon taken from \cite{Koch}. The rate functions are
\begin{equation*}\tag{P1}
\left.\begin{array}{lll}
        a_n(y)=\frac{10-y}{100(\e^{(10-y)/10}-1)}, &\quad a_m(y)=\frac{25-y}{10(\e^{(25-y)/10}-1)}, &\quad a_h(y)= 0.07\e^{-y/20}, \\[2ex]
        b_n(y)=0.125\e^{-y/80}, &\quad b_m(y)=4\e^{-y/18}, &\quad b_h(y)=\frac{1}{\e^{(30-y)/10}+1},
       \end{array}\right.
\end{equation*}
and the constants used are
\begin{equation*}\tag{P1}
\begin{array}{lll}
 E_{\textnormal{Na}}=115 \textnormal{ mV},&\quad \nbar g_{\textnormal{Na}}= 4\textnormal{ pS},&\quad \eta_{\textnormal{Na}}=300\ \mu\textnormal{m}^{-2},\\
 E_{\textnormal{K}}=-12 \textnormal{ mV},&\quad \nbar g_{\textnormal{K}}= 18\textnormal{ pS},&\quad \eta_{\textnormal{K}}=30\ \mu\textnormal{m}^{-2},\\
 E_\textnormal{L}=0 \textnormal{ mV}, &\quad \nbar g_\textnormal{L}= 0.3 \textnormal{ mS/cm$^2$}, & \\
 C= 1\ \mu\textnormal{F/cm}^2\,. & &
\end{array}\phantom{xxxxx[P1]}
\end{equation*}
where $\eta_\textnormal{Na}$ and $\eta_\textnormal{K}$ denote the $Na^+$ and $K^+$ channel density, that is for the numerical example we considered a model with $300$ $Na^+$--channels and $30$ $K^+$--channels. The input we consider is a monophasic current starting at $t=1$ ms and lasting for $1$ ms of strength $30$ pA, i.e., $I(t)=30\cdot\mathbb{I}_{(1,2]}(t)$.

The rate functions and parameters for the second example considered taken from \cite{Minoetal} are given by
\begin{equation*}\tag{P2}
\left.\begin{array}{ll}
        a_m(y)=\frac{1.872 (25.41-y)}{\e^{(25.41-y)/6.06}-1}, &\quad a_h(y)= \frac{-0.549 (27.74+y)}{1-\e^{(y+27.74)/9.06}}, \\[2ex]
        b_m(y)=\frac{3.973 (21.001-y)}{1-\e^{(y-21.001)/9.41}}, &\quad b_h(y)=\frac{22.57}{1+\e^{(56-y)/12.5}}\,.
       \end{array}\right.\phantom{xxxxxxxx[P2]}
\end{equation*}
and the constants used are
\begin{equation*}\tag{P2}
\begin{array}{lll}
 E_{\textnormal{Na}}=144 \textnormal{ mV},&\!\!\nbar g_{\textnormal{Na}}= 2.569\cdot 10^{-5}\textnormal{ pS},&\!\!\eta_{\textnormal{Na}}=1000\ \mu\textnormal{m}^{-2},\\
 E_{\textnormal{K}}=0 \textnormal{ mV},&\!\!\nbar g_{\textnormal{K}}= 0\textnormal{ pS},&\!\!\eta_{\textnormal{K}}=0\ \mu\textnormal{m}^{-2},\\
 E_\textnormal{L}=0 \textnormal{ mV}, &\!\!\nbar g_\textnormal{L}= 1/(1953.49\cdot 10^{3}) \textnormal{ mS/cm$^2$}, & \\
 C= 0.0714\cdot 10^{-6}\ \mu\textnormal{F/cm}^2\,. & &
\end{array}
\end{equation*}
Note that this examples considers $Na^{+}$ channels only. As input we consider a monophasic current starting at $t=0.1$ ms and lasting for 0.1 ms of strength $35.1\cdot 10^{-6}$ pA, i.e., $I(t)=35.1\cdot 10^{-6}\cdot\mathbb{I}_{(0.1,0.2]}(t)$. These are the parameters for the neuron considered in \cite{Minoetal}. We note that the input current strength for the experiments is incorrectly reported in \cite{Minoetal} and corrected in \cite{Bruce}.


\bibliographystyle{plain}
\bibliography{mybib}

\end{document}